\documentclass[11pt,oneside]{amsart}
\usepackage[hyphens]{url}
\usepackage[foot]{amsaddr}
\usepackage{geometry}
\usepackage[T1]{fontenc}
\usepackage[english]{babel}
\addto\extrasenglish{
  
}
\usepackage{amssymb}
\usepackage{amsmath}
\usepackage{amsthm}
\usepackage{amsfonts}
\usepackage[shortlabels]{enumitem}
\usepackage{graphicx}
\usepackage{mathabx}
\usepackage{eurosym}
\usepackage{dsfont}
\usepackage{csquotes}
\usepackage{xparse}
\usepackage{accents}
\usepackage{color}
\usepackage{xcolor}
\usepackage{blindtext}
\usepackage{thmtools}
\usepackage{mathrsfs}
\usepackage{stmaryrd}
\usepackage{amssymb}
\usepackage{amsmath}
\usepackage{thm-restate}
\usepackage[
    backend=biber,
    style=alphabetic,
    isbn=false,
    maxnames=50,
    minalphanames=6,
    maxalphanames=10
  ]{biblatex}
  \AtEveryBibitem{
  	\clearfield{urlyear}
  	\clearfield{urlmonth}
  }

\usepackage[hyperindex,colorlinks,breaklinks]{hyperref}
\hypersetup{
hypertexnames=false,
breaklinks=true,
colorlinks = true,
linkcolor = {purple},
urlcolor = {blue},
citecolor = {brown}
}
\usepackage[hyphenbreaks]{breakurl}
\usepackage{doi}
\numberwithin{equation}{section}
\newtheorem{lem}{Lemma}[section]
\newtheorem{lemap}{Lemma}

\newtheorem{defi}{Definition}[section]
\newtheorem{theorem}{Theorem}

\newtheorem{thmelit}{Theorem}

\newtheorem{cor}[lem]{Corollary}
\newtheorem{proposition}[lem]{Proposition}
\newtheorem{prop}[lem]{Proposition}

\newtheorem{ass}{A\!\!}

\newtheorem{remark}[lem]{Remark}

\newcommand{\RR}{{\mathbb R}}
\newcommand{\NN}{\mathbb N}
\newcommand{\ZZ}{{\mathbb Z}}

\newcommand{\EE}{\mathbb{E}}
\newcommand{\E}{\mathcal{E}}
\newcommand{\bxi}{{\bar{\xi}}}
\newcommand{\Eg}{\mathscr{E}}
\newcommand{\Xg}{\mathfrak{X}}

\newcommand{\XX}{\mathbb{X}}
\newcommand{\X}{\mathcal{X}}
\DeclareMathOperator{\ext}{\mathbf{Ext}}
\DeclareMathOperator{\nexp}{\mathbf{NonExp}}
\newcommand{\UU}{\mathbb{U}}
\newcommand{\U}{\mathcal{U}}
\newcommand{\be}{{\bar{e}}}

\newcommand{\Y}{\mathcal{Y}}
\newcommand{\TT}{\mathbb{T}}
\newcommand\F{{\mathcal F}}
\newcommand{\s}{\mathbf{s}}

\newcommand*{\N}{\mathcal N}
\renewcommand{\L}{\mathcal{L}}
\newcommand{\M}{\mathcal{M}}
\renewcommand{\P}{\hat{\PP}}
\newcommand{\B}{\mathcal{B}}

\newcommand{\as}{\PP(d\bxi)\text{-a.s.}}

\newcommand{\PP}{{\mathbb P}}
\newcommand{\GG}{\mathbb{G}}

\def\beq{\begin{equation}}
\def\eeq{\end{equation}}
\let\origautoref\autoref
\def\autoref#1{\textbf{\origautoref{#1}}}

\title[A Kesten Stigum theorem for MGWREs with infinitely many types]{A Kesten Stigum theorem for Galton-Watson processes with infinitely many types in a random environment}
\author{Maxime Ligonnière}\address{Institut Denis Poisson UMR 7013, Université de Tours, Université d’Orléans, CNRS France}
\email{maxime.ligonniere@univ-tours.fr}
\thanks{The author received the support of the Chair "Modélisation Mathématique et Biodiversité" of VEOLIA-Ecole Polytechnique-MnHn-FX, and of the ANR projects NOLO (ANR 20-CE40-0015) and RAWABRANCH (ANR-23-CE40-0008), funded by the French ministry of research.}
\address{CMAP, CNRS, INRIA, École polytechnique, Institut Polytechnique de Paris, 91120 Palaiseau, France}
\begin{document}
\begin{abstract}
	{\color{black} In this paper, we study a Galton-Watson process $(Z_n)$ with infinitely many types in a random ergodic environment $\bxi=(\xi_n)_{n\geq 0}$.  We focus on the supercritical regime of the process, where the quenched average of the size of the population grows exponentially fast to infinity. We work under Doeblin-type assumptions coming from \cite{ligonniere_ergodic_2023}, which ensure that the quenched mean semi group of $(Z_n)$ satisfies some ergodicity property and admits a $\bxi$-measurable family of space-time harmonic functions. We use these properties to derive an associated nonnegative martingale $(W_n)$. Under a $L\log(L)^{1+\varepsilon}$-integrabilty assumption on the offspring distribution, we prove that the almost sure limit $W$ of the martingale $(W_n)$ is not degenerate. Assuming some uniform $L^2$-integrability of the offspring distribution, we prove that conditionally on $\{W>0\}$, at a large time $n$, both the size of the population and the distribution of types correspond to those of the quenched mean of the population $\EE[Z_n|\bxi, Z_0]$.
	We finally introduce an example of a process modelling a population with a discrete age structure. In this context, we provide more tractable criterions which guarantee our various assumptions are met.}
\end{abstract}
\maketitle
\setcounter{tocdepth}{1}
\tableofcontents
\newpage
\section{Introduction}
Let $(\E, \Eg)$ be a measurable set and $\bxi=(\xi_n)_{n\in{\ZZ_+}}$ be an ergodic sequence with values in $\E$. We refer to $\E$ as the set of environments and to $(\xi_n)_n$ as the environmental process, and note $\theta:\E^{\ZZ_+}\mapsto \E^{\ZZ_+}$ the shift application. Let $(\XX,\X)$ be another measurable space, such that $\{x\}\in \X$ for any $x\in\XX$. We define
\( \Xg=\bigcup_{d\geq 0} \XX^d,\)
which we endow with the smallest $\sigma$-algebra generated by all the $\X^{\otimes d}$, $d\geq 0$. $\XX$ represents the possible types of the individuals. 
We also define the set of point measures 
\[ \N=\left\{\left. \sum_{i=1}^d \delta_{x_i} \right|d\geq 0, (x_1,\dots, x_n) \in \XX^n\right\},\]
and endow it with the $\sigma$-algebra generated by the maps $\mu\in\N\mapsto \mu(A)\in{\ZZ_+}$, for any $A\in\X$.
A population with individuals of various types can either be represented as a tuple of elements of $\XX$ (i.e. by an element of $\Xg$), or by a point measure, that is an element of $\N$. The map \begin{equation} \psi :(x_1,\dots, x_d)\in\Xg \mapsto \sum_{i=1}^d \delta_{x_i}\in \N \label{eq:psi} \end{equation} creates a correspondance between those two representations.

We note $\M_1(\Xg)$ the set of probability measures on $\Xg$, and consider, for each $(x,e)\in\XX\times\E$, a probability distribution $\L_{x,e}\in\M(\Xg)$, which represents the law of the random offpsring of an individual of type $x$ living in the environment $e$.
We consider a probability space $\Omega$, on which are defined both the random process $(\bxi)$ as well as an array $\left(L_{x,e}^{k,n}: x\in \XX,e\in \E, k,n\geq 0\right)$ of independent random variables with values in $\Xg$, independent $\bxi$, such that $L_{x,e}^{k,n}$ is distributed according to $\L_{x,e},$ for any $(x,e)\in\XX\times\E$. We moreover assume that $(x,e,\omega)\in \XX\times \E\times \Omega \mapsto L_{x,e}^{k,n}(\omega)$ is measurable for all $k,n$. Finally, we note $N_{x,e}^{k,n}=\psi\left( L_{x,e}^{k,n}\right),$ and define for each $N\geq 0$  the $\sigma$-algebra $\F_N$ generated by the variables $(L_{x,e}^{k,n})_{x\in \XX,e\in \E, k\geq 0,0\leq n\leq N-1}$.

\begin{defi}
	We call multitype Galton-Watson process $(Z_n)$ in ergodic environment $\bxi=(\xi_n)_{n\geq 0}$ with initial population $Z_0\in \N$ an $\mathcal{N}$-valued process defined through the induction relation
	\begin{equation}\label{eq:MGWRE}Z_{n+1}=\sum_{x\in\XX} \sum_{k=1}^{Z_n(\{x\})} N_{x,\xi_n}^{k,n}.\end{equation}
\end{defi}
For short, we note such a process MGWRE in the rest of the article.
As often when studying random processes in random environment, we note $\PP$ the annealed law, that is the joint law of $\bxi$ and $(Z_n)$, and $\PP_{\bxi}$ the quenched law of $(Z_n)$, obtained by conditioning $\PP$ on the environmental process $\bxi$. 

We define the quenched first moment operators $M_{k,n}$ by setting  
\begin{equation*} M_{k,n}(f)(x)=\EE[Z_n(f)|\bxi,Z_k=\delta_x]=\EE_\bxi\left[Z_n(f)|Z_k=\delta_x\right],\end{equation*}
for any $x\in\XX$ and any bounded function $f$ on $\XX$, and  $\mu M_{k,n}(f)=\int M_{k,n}(f)(x) \mu(dx)$, for any signed measure $\mu$ on $\XX$. The operators $(M_{k,n})_{k\leq n}$ act therefore both on the space of bounded measurable functions $\B(\XX)$, endowed with the supremum norm $\Vert\cdot\Vert_{\infty}$, and on the space of signed measures $\M(\XX)$, endowed with the total variation norm $\Vert \cdot \Vert_{TV}$.
They satisfy the semi group relation
\[M_{k,N}=M_{k,n}M_{n,N}\]
for any $k\leq n\leq N$, and noting $M_{n}=M_{n,n+1}$, it holds
\[M_{k,n}=M_{k}\cdots M_{n-1}.\]
Note that $(M_n)_{n\geq 0}$ are random operators, since they depend on the random environmental process $\bxi$. More precisely, the random operator $M_{k,n}$ is $\sigma(\xi_k,\dots \xi_{n-1})$ measurable, and $(M_{n})_n$ is an ergodic sequence of operators. We note $m_{k,n}(x)=M_{k,n}(\mathds{1})(x)$ for any $x\in\XX$ and $k\leq n$.
\subsection{Preliminary assumptions}
In the sequel, we work under the following assumptions from the article \cite{ligonniere_ergodic_2023}.
\begin{ass} \label{ass: boundedness_KS}
	For almost all $e\in\E$, the function $x\mapsto m_{0,1}(x,e)$ is a bounded and positive function.
\end{ass}
By stationarity, \autoref{ass: boundedness_KS} implies that $\PP$-almost surely, the product $M_{k,n}$ is a continuous, non zero, positive linear operator.

\begin{ass} \label{ass: moments_m reinforced_KS}
	$\EE \left[ \log \Vert m_{0,1} \Vert_\infty \right]=\EE \log \left(\vvvert M_{0,1} \vvvert \right)<\infty,$
\end{ass}
where $\vvvert\cdot\vvvert$ is the operator norm defined as 
\( \vvvert M \vvvert = \sup_{x\in\XX} \sup_{{f\in\B(\XX), \Vert f \Vert_\infty=1}} \vert M(f)(x)\vert = \Vert M\mathds{1}\Vert_\infty.\)
In particular, the submultiplicativity of the norm $\vvvert \cdot \vvvert$, the ergodicity of $(\xi_n)_n$ and \autoref{ass: moments_m reinforced_KS} imply that $\EE \log^+ \left(\left\vvvert M_{k,n} \right\vvvert\right)<\infty$ for all $k\leq n$.
As a consequence of \autoref{ass: moments_m reinforced_KS}, one defines the Lyapunov exponent $\lambda$ of the random sequence of operators $(M_n)_{n\geq 0}$ as
\begin{equation*}
	\lambda=\inf_{n\geq 1}n^{-1}\EE \left[ \log \vvvert M_{0,n}\vvvert \right]\in[-\infty,+\infty)
\end{equation*}
and prove, using submultiplcativity arguments, that almost surely $\lim_{n\rightarrow \infty} \left(\vvvert M_{0,n} \vvvert\right)^{\frac{1}{n}}={e^\lambda}.$

We call admissible coupling constants or admissible triplet a measurable map $(\nu,c,d):\bar{e}=(e_n)_{n\geq 0} \in \mathcal{E}^\NN\mapsto (\nu_{\bar{e}},c({\bar{e}}), d({\bar{e}}))\in\mathcal{M}_1(\XX)\times[0,1]^2 $ such that, $\PP(d\bxi)$-a.s.
\begin{itemize}
	\item[i)]for all $x\in\XX$ and all $f\in \B_+(\XX)$, the couple $(\nu_{\bxi},c(\bxi))$ satisfies \begin{equation}\label{eq:prop c_KS} \delta_x M_{0,1}(f)\geq c(\bxi) \Vert \delta_x M_{0,1} \Vert \nu_\bxi(f)\end{equation}
	\item[ii)] for all $n\geq 0$, the couple $(\nu_\bxi,d(\bxi))$ satisfies 
	\begin{equation} \nu_\bxi (m_{1,n}) \geq d(\bxi) \vvvert M_{1,n} \vvvert \label{eq: prop d_KS}\end{equation}
\end{itemize}
When an admissible triplet is defined, we define the random variable \[ \gamma=c(\bxi)d(\bxi).\]
Note that taking $c=d=0$ and any measurable map $\bar{e}\mapsto \nu_{\bar{e}}$ defines an admissible triplet, however in this case $\gamma=0$. 
Our main assumption, which prevents this trivial situation, is therefore
\begin{ass}
	\label{ass : moments_cd_KS} There exists a triplet of admissible coupling constants $\bar{e}\mapsto(\nu_{\bar{e}},c(\bar{e}),d(\bar{e}))$ such that $\EE[-\log(\gamma)]<\infty$.
\end{ass}
We introduce 
\[\tilde{\eta}=\exp(\EE\left[\log(1-\gamma)\right])\in[0,1)\]
and notice that under \autoref{ass : moments_cd_KS}, we have $0\leq \tilde{\eta}<1$.
As a consequence of results from \cite{ligonniere_ergodic_2023}, assumptions \autoref{ass: boundedness_KS} to \autoref{ass : moments_cd_KS} guarantee that the semi-group $(M_{k,n})$ satisfies some geometric ergodicity properties at rate $\tilde{\eta}$ and a law of large numbers in the form $\lim_{n\rightarrow \infty} \left(\Vert \mu M_{0,n}\Vert_{TV}\right)^{\frac{1}{n}} ={e^\lambda}$.
Moreover, these assumptions guarantee the almost sure existence of 
\[h_k(x)=\lim_{n\rightarrow \infty}\frac{m_{k,n}(x)}{\Vert m_{k,n}\Vert}\in(0,1]\]
for all $x\in\XX$.
We provide in Section \ref{subs:estimate semigroup}, Theorem \ref{thme: approx_KS} a more complete and precise statement of the results from \cite{ligonniere_ergodic_2023} that we use here.
\\ The Lyapunov exponent ${e^\lambda}$  governs the asymptotic quenched size of the population. It can be used to extend the classification of Galton-Watson processes into three regimes to the case of MGWREs. When ${e^\lambda}<1$ (in which case we say that the process is subcritical), it holds $\lim_{n\rightarrow +\infty} \Vert \mu M_{0,n}\Vert=0$ a.s. for any measure $\mu$. Under some additional assumptions, when ${e^\lambda}=1$ (that is, in the critical regime), it is proven in \cite{ligonniere_ergodic_2023} that $\liminf_{n\rightarrow +\infty} \Vert \mu M_{0,n}\Vert=0, \PP(d\bxi)$-a.s.  A classical first moment argument guarantees that the property $\liminf_{n\rightarrow +\infty} \Vert \mu M_{0,n}\Vert=0, \PP(d\bxi)$-a.s is sufficient for the extinction event $\ext=\{\exists n\geq 0, Z_n=0\}$ to have probability $1$. Therefore extinction is a.s. in the subcritical and the critical regime for an MGWRE. In our paper, we rather focus on the supercritical case and assume
\begin{ass}[Supercriticality] \label{ass:supercritical}
	$ \lambda>0$.
\end{ass}
We now present our results.

\subsection{Results}
In the supercritical case $\lambda>0$, one expects that the population survives on an event of positive probability. However, when studying MGWREs with infinitely many types this is not as easy to prove as when $\XX$ is finite, and is, to the best of our knowledge, not proven in full generality. In this article, we obtain in particular, under some uniform moment assumption, the survival of the process in the supercritical case, as well as a description of the surviving population, in the form of on a Kesten-Stigum type theorem. The original Kesten-Stigum result was proven in \cite{kesten_galton-watson_1966} and deals with Multitype Galton-Watson processes with finitely many types, in a constant environment, an alternative proof based on a classical change of measure which inspires some of our techniques is presented in \cite{lyons_conceptual_1995}. Many similar results have been developed since on various random models. For example, in the case of general time homogeneous branching processes with infinitely many types, similar results appear in \cite{asmussen_strong_1976,athreya_change_2000}. \cite{englander_strong_2010, louidor_strong_2020, medous_spinal_2024} also present related results on some branching models, respectively branching diffusions, branching brownian motion with absorbtion, and some continuous space branching process with interactions. Closer to our context, Kesten-Stigum theorems on multitype Galton-Watson processes in random environment with finitely many types have recently been established in \cite{grama_kestenstigum_2023,grama_limit_2024}. Roughly speaking, these results deal with four problems :
\begin{itemize}[-]
	\item Exhibiting a martingale of the form 
	$W_n=Z_n(f_n)/r_{0,n}$
	where $f_n\in \mathcal{B}_+(\XX)$ and $r_{0,n}\in \RR_+$ are $\bxi$-measurable (deterministic in constant environment)
	\item Under some $L\log(L)$-integrability condition on the offspring distribution, proving that the almost sure limit $W:=\lim_{n\rightarrow \infty} W_n$ is not degenerate, in the sense that $\PP[W>0]>0$ (in random environment, this statement typically becomes $\PP_{\bxi}(W>0)>0, \PP(d\bxi)$-a.s.).
	\item Proving that the extinction event $\bigcup_{n\geq 0} \{Z_n=0\}$ coincides with the event $\{W=0\}$
	\item Obtaining some almost sure asymptotic estimate on $Z_n$ describing the distribution of the types in the population conditionally on the event $\{W>0\}.$
\end{itemize}
We adress those four topics in this article. Our respective contributions regarding each of these points are described in following four subsubsections.
\subsubsection{The fundamental martingale}
Our first result generalizes the martingale obtained in \cite{grama_kestenstigum_2023} to the case of an infinite type set $\XX$.
For all $k\geq 0$, we note $\lambda_k=\Vert M_k h_{k+1}\Vert_\infty$ and $\lambda_{k,n}=\lambda_k\dots \lambda_{n-1}$ for all $k\leq n$.
Moreover, we set

$$W_n=\frac{Z_n(h_n)}{\lambda_{0,n}Z_0(h_0)}.$$
Our first main result is
\begin{theorem}[The fundamental martingale]\label{thme: martingale}
	Consider a MGWRE satisfying assumptions \autoref{ass: boundedness_KS} to \autoref{ass : moments_cd_KS}. Then, conditionally on the environmental process $(\xi_n)_{n\geq 0}$, the process $(W_n)$ is a $(\F_n)_{n\geq 0}$-martingale.
\end{theorem}
The main argument behind this theorem is the space-time harmonicity of the family $(h_k)$, in the sense that $M_{k}h_{k+1}=\lambda_k h_k$ for all $k$. This is the core of Lemma \ref{lem:eigenfunction}, stated in Subsection \ref{subs:proof_martingale}.
\\When $\XX$ is finite, the counterpart of $h_k$ is the limit $v_k$ of the sequence of right dominant eigenvectors $(v_{k,n})_{n\geq k}$ of the matrices $(M_{k,n})_{n\geq k}$.  As a consequence, our martingale $W_n$ is similar to the one obtained in finite dimension by \cite{grama_kestenstigum_2023}.
\\ Since the sequence $(\lambda_n)$ is ergodic, it is straightforward from Birkhoff's ergodic theorem that $(\lambda_{0,n})^\frac{1}{n}$ converges almost surely to a constant. We additionally prove in Subsection \ref{subs:proof_martingale}, Prop \ref{prop: prod lambda} that this limit is actually ${e^\lambda}$.

\subsubsection{Non degeneracy of the martingale}
By Theorem \ref{thme: martingale}, the process $(W_n)$ is a positive martingale, therefore there exists a random variable $W\in[0,\infty)$ such that it holds
\[W_n\underset{n\rightarrow \infty}{\longrightarrow} W \quad \PP-\text{a.s.}\]
We wonder now whether $W=0$ almost surely or $\{W>0\}$ with positive probability. To perform this study of $W$, we introduce the assumption
\begin{ass} \label{ass: LlogL} There exists a nondecreasing and positive function $f$ such that $\frac{1}{x\log(x) f(x)}$ is integrable at $+\infty$ and \[ \EE \left[ \sup_{x\in \XX} \frac{ \EE_{\bxi}\left[ Z_1(h_1) \log[ Z_1(h_1)] f(Z_1(h_1))| Z_0=\delta_x \right] }{M_0h_1(x)}\right] <\infty.\]\end{ass}
Note that the functions $f(t)=\log^+(t)^\varepsilon$, for any $\varepsilon>0$ satisfy the integrability assumption $\int^{+\infty} \frac{dx}{x\log(x) f(x)}<+\infty$.
Let us now state our main result dealing with the limit of the martingale.
\begin{theorem}[Non-degeneracy of the fundamental martingale]\label{thme:nondegenerescence}
	Consider a MGWRE satisfying assumptions \autoref{ass: boundedness_KS} to \autoref{ass: LlogL}.
	For any initial population $Z_0\in \N$, it holds
	\[\PP_\bxi[W>0|Z_0]>0\text{ and } \EE_\bxi[W|Z_0]=1 \quad \PP(d\bxi)\text{-a.s.}\]	
\end{theorem}
Since $W=0$ on the extinction event, our result implies in particular that the population survives with positive probability. 
The non degeneracy of $(W_n)$ usually requires a $L\log L$ moment on the offspring distribution, see \cite{kesten_galton-watson_1966,lyons_conceptual_1995} for example. Our assumption \autoref{ass: LlogL} adapts this condition to our context. Indeed, setting $L=Z_1(h_1)$, \autoref{ass: LlogL} can be rewritten 
\[ \EE\left[ \sup_{x\in \XX} \frac{\EE_\bxi[ L \log^+(L)f(L)|Z_0=\delta_x]}{\EE_\bxi \left[ L |Z_0=\delta_x\right]}\right]<\infty.\]
\begin{remark} When $\XX$ is finite, the supremum over the type space $\XX$ can be dominated by a sum over all the elements of $\XX$, thus condition \autoref{ass: LlogL} reduces further to
	\begin{equation} \ \EE\left[\left.\frac{L\log^+(L)f(L)}{\lambda_0 h_0(x)}\right|Z_0=\delta_x\right]<\infty \label{eq: Llog(L) X fini} \end{equation}
	for any $x\in\XX$.
	This condition is very close to condition H6 of \cite{grama_kestenstigum_2023}. Indeed, this condition H6 is also an integrability condition on the quantity $L=Z_1(h_1)$, which consists in giving a weigh $h_1(x)$ to individuals of type $x$ alive at time $1$. The only difference between \eqref{eq: Llog(L) X fini} and condition H6 in \cite{grama_kestenstigum_2023} is  the multiplicative term $f(L)$, which is only a very slight strengthening of the condition since one can choose for example $f(t)=\log^+(t)^\varepsilon$.
\end{remark}
In the case of an infinite type set $\XX$, the presence of the supremum over $x\in\XX$ is the main difficulty in Assumption \autoref{ass: LlogL}. We introduce some examples of offspring distributions with various decay rates which satisfy \autoref{ass: LlogL} in Subsection \ref{subs: proof_nondegeneracy}.
\subsubsection{Asymptotic type distribution}
Our next result deals with the asymptotic type distribution in the population, in the form of almost sure asymptotic estimates of $Z_n(f)$, for a large classe of functions $f$. It relies on the following uniform second moment assumption :
\begin{ass} \label{ass:L2}
	\[\EE\left[\log^+\left(\sup_{x\in\XX} \EE_\bxi \left[\left.\Vert Z_1\Vert^2\right| Z_0=\delta_x \right]\right)\right]<\infty.\]
\end{ass}
Before stating the theorem, we introduce the notation $\log^-(x)= -\min(\log(x),0)>0$, so that $\log(x)=\log^+(x)-\log^-(x)$.
\begin{theorem}\label{thme:type_distrib} Assume \autoref{ass: boundedness_KS} to \autoref{ass:L2}. Then 
	\begin{enumerate}[i)]
		\item For any $f\in\mathcal{B}(\XX)$, for any $\varepsilon >0$, it holds
		\begin{equation} \label{eq: asymptotic development Z_n(f)}Z_n(f)= W_n  Z_0 M_{0,n}(f) + \underset{n\rightarrow +\infty}{o}\left(\left(\max\left(\tilde{\eta}{e^\lambda},e^{\frac{\lambda}{2}}\right)+\varepsilon\right)^n\right) \quad \PP\text{-a.s.}\end{equation}
		\color{black}	\item For any $f\in\mathcal{B}(\XX)$, $f\geq 0$, if $\EE\left[\log^-\left(\nu M_{1,N} f\right) \right]<+\infty$ for some $N\geq 1$ then 
		\begin{equation} \label{eq: asymptotic development Z_n(f) 2} \lim_{n\rightarrow \infty}\frac{Z_n(f)}{ Z_0 M_{0,n}(f)}= W \quad \PP\text{-a.s}\end{equation}
		and moreover, conditionally on $\{W>0\}$,
		$$\lim_{n\rightarrow \infty} \left( Z_n(f) \right)^\frac{1}{n}={e^\lambda}.\quad \PP\text{-a.s.}$$
	\end{enumerate}
\end{theorem}
\color{black}
Note that \eqref{eq: asymptotic development Z_n(f) 2} only derives from \eqref{eq: asymptotic development Z_n(f)} for functions $f$ such that \begin{equation}\label{eq: condition f}\liminf_{n\rightarrow \infty} \left(Z_0 M_{0,n}(f)\right)^{\frac{1}{n}}>\max\left(\tilde{\eta}{e^\lambda}, e^{\frac{\lambda}{2}} \right).\end{equation} \eqref{eq: condition f} is a straightforward consequence of the properties of the Lyapunov exponnent ${\lambda}$ if $f>c>0$ for some constant $c$, however we do not see a reason for it to hold for a general bounded function $f$.
In Subsection \ref{subs: Leslie asymptotic type}, we focus on the example of a GWRE with countably many types representing a population structured in age and prove that under some mild integrability assumption on the mean matrix, any nonnegative and nonzero function $f$ actually satisfies the assumption $\EE\left[\log^-\left(\nu M_{1,N} f\right) \right]<+\infty$ and thus \eqref{eq: asymptotic development Z_n(f) 2}.
\color{black}
\\ Theorem \ref{thme:type_distrib} allows us to describe both the size and the distribution of types within the surviving population on the event $\{W>0\}$, as stated in the following corollary. 
\begin{cor}\label{cor:size and types} Assume \autoref{ass: boundedness_KS} to \autoref{ass:L2}. Then $\PP$-a.s., conditionally on $\{W>0\}$, it holds
	\begin{equation} \label{eq: size pop} Z_n(\mathds{1})\underset{n\rightarrow \infty}{\sim} W_n \Vert Z_0 M_{0,n} \Vert_{TV} \text{ and } \lim_{n\rightarrow \infty} \left( Z_n(\mathds{1})\right)^\frac{1}{n} ={e^\lambda}.\end{equation}
	Additionally, for any $f\in \mathcal{B}(\XX)$, conditionally on $\{W>0\}$, it holds
	\begin{equation} \label{eq: freq types pop} \lim_{n\rightarrow \infty} \frac{Z_n(f)}{Z_n(\mathds{1})}-\frac{Z_0 M_{0,n}(f)}{Z_0 M_{0,n}(\mathds{1})}=0 \quad \PP\text{-a.s.}\end{equation}
\end{cor}
Claim \eqref{eq: size pop} is obtained by plugging $f=\mathds{1}$ in the previous Theorem. It proves that the population sizes grows almost surely at rate ${e^\lambda}$ on the event $\{W>0\}$. Moreover, applying the Theorem to various test functions, we obtain \eqref{eq: freq types pop}, which provides insight on the asymptotic distribution of types in the population. As $n$ is large, on the event $\{W>0\}$, the respective frequencies of any given type in the population $Z_n$ and in the quenched mean $\EE_\bxi[Z_n|Z_0]$ become similar. 

\color{black}

\subsubsection{Extinction and explosion}
It is common when studying the non degeneracy of the martingale $(W_n)$ of a Galton-Watson process to prove additionally that the non-explosion event $\nexp=\{W=0\}$ coincides with the extinction event $\ext=\{\lim \Vert Z_n \Vert =0\}$ up to a negligible event. In the single-type case, or when studying a multitype GW process with a finite number of types, in constant environment, this relies on noticing  that the probability of extinction $\PP[\ext]$ and the probability of non-explosion $\PP[\nexp]$ both are a fixed point of the generating function of the reproduction law. A convexity argument allows to show then that $\PP[\ext]=\PP[\nexp].$ The obvious inclusion $\ext\subset \nexp$ concludes the proof. This method can be adapted in ergodic environment, with a variant on the fixed point argument, see \cite[p.27]{grama_kestenstigum_2023} and \cite[Prop 3.1]{kaplan_results_1974}. However, when considering an infinite type set, and even in fixed environment, the generating function might have several fixed points, see e.g. \cite{braunsteins_extinction_2018,braunsteins_extinction_2019,braunsteins_pathwise_2019}. Therefore we do not expect that $\ext$ and $\nexp$ coincide in general when $\XX$ is infinite.
We are however able to prove that they do coincide when the population of a specific type $x_0$ is unbounded a.s. on survival.
\begin{prop} \label{prop:extinction-explosion}
	Let $(Z_n)$ be a MGWRE. We assume that \autoref{ass: boundedness_KS}-\autoref{ass: LlogL} hold, that the environmental sequence $\bxi$ is i.i.d and that there exists a type $x_0\in\XX$ such that $\PP$-a.s, conditionally on the survival event, for any initial population $Z_0$,
	\begin{equation}\label{eq:survie type} \limsup_{n\rightarrow \infty} Z_n(\{x_0\})=+\infty.\end{equation}
	Then, for any initial population $Z_0$, the events $\ext=\{\exists n\geq 0, Z_n=0\}$ and $\{W=0\}$ coincide up to a $\PP$-negligible event.
\end{prop}
We provide in subsection \ref{subs: Leslie extinction-explosion} a specific example where assumption \eqref{eq:survie type} holds.
\subsection{Structure of the paper}
In section \ref{sec:prelim}, we introduce some preliminary results on the quenched mean of the process from a previous article, which are useful in the rest of the paper. 
Theorem \ref{thme: martingale} is proved in subsection \ref{subs:proof_martingale}, Theorem \ref{thme:nondegenerescence} is proved in subsection \ref{subs: proof_nondegeneracy}, Theorem \ref{thme:type_distrib} is proved in subsection \ref{subs:proof_type distrib} and Proposition \ref{prop:extinction-explosion} is proved in subsection \ref{subs:extinction-explosion}. Finally, in Section \ref{sec:LGWRE}, we introduce an example of a GWRE with countably many types, representing a population structured in age. This example is an interesting toy model on which we apply our result and discuss our assumptions, in particular Assumption \autoref{ass: LlogL} and criterion \eqref{eq:survie type}.
\section{Preliminaries : behavior of the quenched mean of the process}
\label{sec:prelim}
{\color{black}As often when studying Galton-Watson processes in random environment, a prerequisite is a fine understanding of the quenched mean of the process. In our case, this quenched mean is the random semi group of nonnegative operators $(M_{k,n})$, which act on $\mathcal{M}(\XX)$ on the left and $\mathcal{B}(\XX)$ on the right.} We use the assumptions and notations of \cite{ligonniere_ergodic_2023}. In particular, we recall that $(\nu,c,d)$ is an admissible triplet, that $\gamma=c(\bxi)d(\bxi)$ and that $\theta$ refers to the shift mapping on $\mathcal{E}^{\NN}$. We introduce the random variables
\[c_n=c(\theta^n(\bxi)), d_{n+1}=d(\theta^n(\bxi)), \gamma_n=c_nd_{n+1}, \nu_n=\nu_{\theta^n(\bxi)}.\]
in such a way that
\[\nu_0=\nu, c_0=c, d_1=d, \gamma_0=\gamma.\]
The main result we rely on in this paper, which was proven in \cite{ligonniere_ergodic_2023}, i.e. \cite{ligonniere_ergodic_2023}, is the following.
\label{subs:estimate semigroup}
\begin{thmelit}[\cite{ligonniere_ergodic_2023}] \label{thme: approx_KS}
	Under assumptions \autoref{ass: boundedness_KS} to \autoref{ass : moments_cd_KS}, the semi group $(M_{k,n})_{k\leq n}$ satisfies, $\PP(d\bxi)--$a.s.,
	\begin{enumerate}[i)]
		\item For any $k\geq 0$, the sequence $\left(\frac{m_{k,n}}{\Vert m_{k,n} \Vert}\right)_{n\geq k}$ converges uniformly on $\XX$ to some random limit $h_k\in\mathcal{B}(\XX)$.
		\item For any $\eta>\tilde{\eta}$, for any $n\geq k\geq 0$ and any $\mu_1,\mu_2 \in \M_+(\XX)$, 
		\begin{equation}\label{eq:PF1_KS}
			\left \Vert \mu_1 M_{k,n} - \frac{\mu_1(h_k)}{\mu_2(h_k)} \mu_2 M_{k,n} \right\Vert_{TV}\leq \frac{4}{\gamma_{n-1}} \prod_{i=k}^{n-1} (1-\gamma_i) \Vert \mu_1 M_{k,n}\Vert=\underset{n\rightarrow \infty}O(\eta^n\Vert \mu_1 M_{k,n}\Vert).
		\end{equation}
		\item For any finite, positive and non-zero measure $\mu$ on $\XX$,
		\begin{equation} \label{eq : growth rate_KS}
			\underset{n\rightarrow \infty}\lim \frac{1}{n}\log \Vert \mu M_{0,n} \Vert =\underset{n\rightarrow \infty}\lim\frac{1}{n} \log \vvvert M_{0,n}\vvvert =\inf_{n\geq 1}\frac{1}{n}\EE \left[ \log \vvvert M_{0,n}\vvvert \right]=\lambda\in[-\infty,\infty).
		\end{equation}
	\end{enumerate}
\end{thmelit}
This provides a notion of ergodicity for the first moment semi-group $(M_{k,n})_{k\leq n}$, in the sense that the measure $\mu_1 M_{k,n}$ asymptotically only depends on $\mu_1$ through a multiplicative factor $\mu_1(h_k)$. {\color{black}It generalizes both Perron-Frobenius' theorem and previous results of \cite{hennion_limit_1997} dealing with products of random positive 
	$d\times d$ matrices to the case of products of random infinite dimensional operators. }
In particular, by \textit{ii)}, choosing an arbitrary measure $\pi_0\in \mathcal{M}_+(\XX)$ and setting $\pi_n=\frac{\pi_0 M_{0,n}}{\Vert \pi_0 M_{0,n}\Vert}$ defines a sequence of random probability measures such that $\PP(d\bxi)$-a.s., for any $\mu\in \mathcal{M}_+(\XX)$ and any $\eta\in (\tilde{\eta},1)$, it holds
\begin{align}\left\Vert \mu M_{k,n} - \frac{\mu(h_k)}{\pi_0 M_{0,k} h_k} \Vert \pi_0 M_{0,n} \Vert_{TV} \pi_n \right \Vert_{TV} &\leq \frac{4}{\gamma_{n-1}} \prod_{i=k}^{n-1}(1-\gamma_i)\Vert \mu M_{k,n} \Vert\label{eq: PF2}\\ &=\underset{n\rightarrow\infty}{O}(\eta^n\Vert \pi_0 M_{0,n} \Vert_{TV})\label{eq: PF2}.\nonumber\end{align}
In the context of multitype Galton-Watson processes, Theorem \ref{thme: approx_KS} thus implies that, as $n\rightarrow \infty$, \begin{itemize}[-]
	\item the size of the quenched mean of the population $\Vert \mu M_{k,n}\Vert_{TV} = \left \Vert\EE[Z_n|Z_k=\mu, \bxi] \right\Vert_{TV}$ is close to $\Vert \pi_0 M_{0,n} \Vert_{TV}$, up to a multiplicative factor which is constant in time and depends linearly on $\mu$,
	\item the quantity $\Vert \pi_0 M_{0,n} \Vert_{TV}$ is independent of $\mu$ and is of order $e^{n\lambda}$ for large $n$, up to a subexponential factor,
	\item  the direction $\frac{\mu M_{k,n}}{\Vert \mu M_{k,n}\Vert_{TV}}$ of the measure $\mu M_{k,n}$, which represents the distribution of the types in the mean population at time $n$, is close to the probability $\pi_n$. Asymptotically, it does therefore not depend on $\mu$.
\end{itemize}
As a consequence, for any initial measure $Z_0$, any fonction $f\in\mathcal{B}(\XX)$ and any $\varepsilon>0$, it holds $\PP(d\bxi)$-a.s.
\[Z_0 M_{0,n}(f)=\frac{Z_0(h_0)}{\pi_0(h_0)} \pi_0 M_{0,n}(f)+o\left((\tilde{\eta}{e^\lambda}+\varepsilon)^n\right).\]
In particular, this yields the following corollary of Theorem \ref{thme:type_distrib}.
\begin{cor}
	Under assumptions \autoref{ass: boundedness_KS} to \autoref{ass:L2}, for any $f\in\mathcal{B}(\XX)$, it holds for any $\varepsilon>0$,
	$$Z_n(f)=W_n \frac{Z_0(h_0)}{\pi_0(h_0)} \pi_0 M_{0,n}(f) +o\left(\left(\max(\tilde{\eta}{e^\lambda},e^{\frac{\lambda}{2}})+\varepsilon\right)^n\right) \quad \PP\text{-a.s.}$$
\end{cor}
In \cite{ligonniere_ergodic_2023}, it is additionally proven that the random sequence of probability measures  $\left(\frac{\pi_0 M_{0,n}}{\Vert \pi_0 M_{0,n}\Vert_{TV}}\right)_{n\geq 0}$ converges weakly to a random measure $\Lambda$ on the space $\mathcal{M}_1(\XX)$.

\section{Proofs}
\label{sec: proofs}
\subsection{The fundamental martingale}
\label{subs:proof_martingale}
The main reason for which $(W_n)$ is a martingale is the fact that the sequence $(h_k)$ constitutes a family of space time harmonic functions, in the sense of \cite{biggins_multi-type_1999}.
\begin{lem} \label{lem:eigenfunction}
	It holds $M_{k}h_{k+1}=\lambda_k h_{k}$ $\PP(d\bxi)$-a.s. for all $k\geq 0$. \end{lem}
\begin{proof}
	Let us fix $k\geq 1$. Applying \eqref{eq:PF1_KS} to the constant function $\mathds{1}$, one proves that almost surely, for any measures $\mu_1,\mu_2$, it holds
	\[ \mu_1(m_{k,n}) \underset{n\rightarrow \infty}\sim \frac{\mu_1(h_k)}{\mu_2(h_k)} \mu_2 (m_{k,n}).\]
	Therefore, 
	\[\frac{\mu_1(h_k)}{\mu_2(h_k)}=\underset{n\rightarrow\infty}{\lim} \frac{\mu_1(m_{k,n})}{\mu_2(h_{k,n})}.\]
	As a consequence, applying the previous results to $\mu_1 M_{k-1}$ and $\mu_2 M_{k-1}$, we obtain
	\[\frac{\mu_1 M_{k-1} h_k}{\mu_2 M_{k-1}h_k}=\underset{n\rightarrow\infty}{\lim} \frac{\mu_1(m_{k-1,n})}{\mu_2(h_{k-1,n})}=\frac{\mu_1(h_{k-1})}{\mu_2(h_{k-1})}.\]
	Thus with $\mu_1=\delta_x$ and an arbitrary $\mu_2\neq 0$, we derive
	\[  M_{k-1} h_k(x)= h_{k-1}(x)  \frac{\mu_2 M_{k-1}h_k}{\mu_2(h_{k-1})}.\]
	The functions $M_{k-1} h_k$ and $h_{k-1}$ are therefore $\PP(d\bxi)$-a.s. colinear. Since $\Vert h_k \Vert_\infty=1$, letting $\lambda_{k-1}=\Vert M_{k-1}h_{k} \Vert_\infty\in(0,\infty)$ yields $ M_{k-1}h_{k}=\lambda_{k-1}h_{k-1}$ as expected.
\end{proof}

\begin{proof}[Proof of Theorem \ref{thme: martingale}]
	To prove that $(W_n)$ is a martingale, we must compute
	\[\EE\left[W_{n+1}| (\xi_k)_{k\geq 0}, \F_n\right]=\EE\left[\left.\frac{Z_{n+1}(h_{n+1})}{\lambda_{0,n+1}}\right| \bxi, Z_n\right].\]
	Notice that $h_{n+1}$ and the coefficients $\lambda_0,\cdots,\lambda_n$ are measurable with respect to the environmental process $(\xi_k)$. Then 
	\begin{align}\EE\left[W_{n+1}| \bxi, \F_n\right]&=\frac{\EE\left[ \left. Z_{n+1}(h_{n+1})\right |\bxi,Z_n \right] }{\lambda_{0,n+1}Z_0(h_0)} \nonumber \\
		&=\frac{\sum_{x\in\XX} \sum_{k=1}^{Z_n(\mathds{1}_x)} \EE\left[ \left. N_{x,\xi_n}^{k,n}(h_{n+1}) \right |\bxi,Z_n \right] }{\lambda_{0,n+1}Z_0(h_0)} \nonumber \\ &=\frac{\sum_{x\in\XX} \sum_{k=1}^{Z_n(\mathds{1}_x)} M_n h_{n+1}(x)}{\lambda_{0,n+1}Z_0(h_0)}.\nonumber \\
		&= \frac{ Z_n M_n h_{n+1}}{\lambda_{0,{n+1}} Z_0(h_0)}\end{align}
	By Lemma \ref{lem:eigenfunction}, this yields
	\[\EE\left[W_{n+1}| \bxi, \F_n\right]=\frac{\lambda_n Z_n(h_n)}{\lambda_{0,n+1}Z_0(h_0)}=\frac{Z_n(h_n)}{\lambda_{0,n} Z_0(h_0)}=W_n.\]
\end{proof}
We provide additionally the following estimate
\begin{prop} \label{prop: prod lambda}
	Under assumptions \autoref{ass: boundedness_KS} to \autoref{ass : moments_cd_KS}, it holds $\lim_{n\rightarrow \infty} \left(\lambda_{k,n}\right)^\frac{1}{n}={e^\lambda}\quad \as$ for any $k\geq 0$.
\end{prop}
\begin{proof}
	Iterating the relation $M_{k}h_{k+1}=\lambda_k h_k$, proven in Lemma \ref{lem:eigenfunction}, we get for any $k\leq n$
	\[M_{k,n}h_n=\lambda_{k,n}h_k.\]
	Thus $\lambda_{k,n}=\vvvert M_{k,n} h_n \vvvert.$ Since $\Vert h_n \Vert_\infty =1$, this proves that $\lambda_{k,n}\leq \vvvert M_{k,n} \vvvert.$
	By \textit{i)} of Theorem \ref{thme: approx_KS}, the function $h_n$ is the uniform limit of the sequence $\left(\frac{m_{n,n+p}}{\Vert m_{n,n+p}\Vert_\infty}\right)$, as $p\rightarrow \infty$. Therefore 
	\[ \lambda_{k,n}=\Vert M_{k,n} h_n \Vert_\infty= \left \Vert M_{k,n} \underset{p\rightarrow \infty}{\lim} \frac{m_{n,n+p}}{\Vert m_{n,n+p} \Vert}\right\Vert=\underset{p\rightarrow \infty}{\lim} \left \Vert \frac{M_{k,n} m_{n,n+p}}{\vvvert M_{n,n+p}\vvvert}\right \Vert_{\infty}=\underset{p\rightarrow \infty}\lim  \frac{\vvvert M_{k,n+p} \vvvert}{\vvvert M_{n,n+p}\vvvert }.\]
	Additionally, if $0\leq k\leq n-1$ then for $p\geq 0$,
	\[ M_{k,n+p}=M_{k,n-1} M_{n-1} M_{n,n+p} \geq m_{k,n} c_{n-1} \nu_{n-1} M_{n,n+p}.\]
	Applying this to the function $\mathds{1}$, and using the definition of $d_n$, we get 
	\[ m_{k,n+p}\geq m_{k,n} c_{n-1} \nu_{n-1}(m_{n,n+p}) \geq m_{k,n} c_{n-1} d_n  \nu_{n-1} \vvvert M_{n,n+p}\vvvert.\]
	Taking now the supremum norm yields
	\begin{equation} \vvvert M_{k,n+p} \vvvert = \Vert m_{k,n+p}\Vert_\infty  \geq \gamma_{n-1} \vvvert M_{k,n} \vvvert \vvvert M_{n,n+p} \vvvert.\label{eq:coef gamma lower bound mult norm}\end{equation}
	Thus
	\[ \lambda_{k,n} = \underset{p\rightarrow \infty }{\lim} \frac{\vvvert M_{k,n+p} \vvvert}{\vvvert M_{n,n+p}\vvvert } \geq \gamma_{n-1} \vvvert M_{k,n} \vvvert.\]
	
	We therefore have, for $k<n$
	\begin{equation} \gamma_{n-1} \vvvert M_{k,n} \vvvert \leq \lambda_{k,n} \leq \vvvert M_{k,n} \vvvert \end{equation}
	The result then derives from Theorem \ref{thme: approx_KS} and \cite[Lemma 3.3]{ligonniere_ergodic_2023}.
\end{proof}
\subsection{Condition for nondegeneracy of the fundamental martingale}
\label{subs: proof_nondegeneracy}
\subsubsection{Construction of the Galton-Watson tree}
In the introduction, we defined the MGWRE as a random measure-valued process. Our proof of the nondegeneracy of $(W_n)$ requires however to consider an underlying tree structure. 
Thus, let us recall first a construction of a multitype Galton-Watson tree in a random environment.
Instead of the array $(L_{x,e}^{k,n})$, we consider a similar array $(L_{x,e}^u)_{x\in \XX, e \in \E, u\in\UU}$ indexed by the Ulam-Harris-Neveu tree $\UU$.  We assume that $L_{x,e}^u \sim \L_{x,e}$ and that the sub-arrays $(L_{x,e}^u)_{e\in\EE}$, for $x\in\XX$ and $u\in\U$ are independent from each other, as well as from $(\xi_n)$. For any $(x,e,u)$, the $\Xg$-valued random variable $L_{x,e}^u$ describes the offspring of an individual $u$ of type $x$ reproducing in the environment $e$. As we did previously, we note $N_{x,e}^u= \psi(L_{x,e}^u)$ its measure valued counterpart. We endow $\UU$ with the $\sigma$-algebra $\U$ of events that only depend on the $n$ first generations of the tree for some $n$, and endow $\XX^\UU$ with the cylindrical $\sigma$-algebra. Then, let us build recursively a random labelled subtree $(\TT,\XX)$, where $\TT\subset \UU$ and $X:\TT\longrightarrow \XX$ as follows.
We note $x\in\XX$ the type of the initial individual and set $X(\emptyset)=x$.  Let us note $\GG_n$ the random set of the individuals of $\TT$ which are at height $n$. Then, conditionally on $\GG_n$ and $X|_{\GG_n}$, we set \[ \GG_{n+1}=\left\{ ui, u\in \GG_n, 1\leq i\leq N_{X(u),\xi_n}^u(\mathds{1}) \right\}.\]
Moreover, conditionally on $N_{X(u),\xi_n}^u(\mathds{1})=d$ we define $(X(u1), \cdots , X(ud))=L^u_{X(u),\xi_n}.$ 

An induction argument shows that the process $(\tilde{Z}_n)$ defined by 
\[  \tilde{Z}_n=\sum_{u \in \GG_n } \delta_{X(u)}\]
is distributed as $(Z_n)$.
Indeed, 
\begin{align*}
	\tilde{Z}_{n+1}&=\sum_{v\in \GG_{n+1}} \delta_{X(v)}= \sum_{u\in \GG_n} \sum_{i=1}^{N_{X(u),\xi_n}^u(\mathds{1})} \delta_{X(ui)}= \sum_{u\in \GG_n} N^u_{X(u),\xi_n}
\end{align*}
Then, using the map $\psi$ defined in \eqref{eq:psi}, we write
\begin{align*}
	\tilde{Z}_{n+1}&= \sum_{x\in\XX} \sum_{\substack{u\in\GG_n\\ X(u)=x}}\psi(L^u_{x,\xi_n}). \end{align*}
For any type $y\in\XX$, conditionally on $\GG_n$ and $\bxi$, the family $A_y=(\psi(L_{x,\xi_n}^u))_{u\in\GG_n, X(u)=y}$ contains $\tilde{Z}_n(\mathds{1}_y)$ independent variables, all of which are distributed according to the $\psi$-pushforward of the distribution $\L_{y,\xi_n}$. Moreover, the families $(A_y)_{y\in\XX}$ are mutually independent conditionally on $\GG_n$ and $\bxi$. Therefore, conditionally on $\tilde{Z}_n$ and $\bxi$, it holds
\[ \tilde{Z}_{n+1}\overset{d}{=} \sum_{x\in\XX} \sum_{k=1}^{\tilde{Z}_n(\mathds{1}_x)} N_{x,\xi_n}^{k,n}.\]
Therefore the law of $\tilde{Z}_{n+1}$ conditionally on $\tilde{Z}_n$ is the same as the law of $Z_{n+1}$ conditionally on $Z_n$, which proves that 
\[ (\tilde{Z}_{n})_{n\geq 0}\overset{d}{=}(Z_n)_{n\geq 0}\]
since $Z_0=\tilde{Z}_0$. In the sequel, by abuse of notation, we shall denote $Z_n=\tilde{Z}_n$.
\subsubsection{Construction of the size biased tree}
We introduce a slight modification in the above construction to create a biased law on the space of labelled trees. More precisely, we build a random labelled tree $(\TT^*, X^*)$, equipped with a random sequence of marked individuals $(\s_n)_{n\geq 0}$, called the spine. 

We recall that for almost any sequence of environments $\be=(e_n)_{n\geq 0}\in \E^{\ZZ_+}$, there exists a sequence of functions $(h_n)\in\B(\XX)^{\ZZ_+}$ and a sequence of positive real numbers $(\lambda_n)$ satisfying Theorem \ref{thme: approx_KS}. We fix such a sequence of environments and introduce, for each type $x$ and almost any sequence $\be$ the biased probability distribution $\L^*_{x,\be}$ on $\Xg$ such that
\begin{equation} \label{eq:biased law} d\L^*_{x,\be}(Z)=\frac{Z(h_1)}{\lambda_0 h_0(x)} d\L_{x,e_0}(Z).\end{equation}
\begin{remark}{\color{black} In case of a finite type set $\XX=\{1,\dots, d\}$, the functions $h_0$ and $h_1$ are respectively the limits of the right dominant eigenvectors of the matrices $M_{0,n}$ and $M_{1,n}$. It holds therefore $\PP[N^*_{x,e}=u]= \frac{\langle u,h_1\rangle}{\lambda_0 h_0(x) } \PP[N_{x,e_0}=u]$ for any vector $u\in \ZZ_+^d$ (with the abuse of notation $\N=\ZZ_+^d$). One recovers the biased distribution used in \cite{grama_kestenstigum_2023}. This distribution is itself an extension of the biased laws introduced in \cite{lyons_conceptual_1995,friedman_conceptual_1997} for Galton Watson processes in constant environment respectively with one and a finite number of types.}
\end{remark}
In particular, for any $f\in \B(\XX)$, it holds 
\[ \int_{\Xg} \psi(Z)(f)d\mathcal{L}^*_{x,e}(Z)=\int \psi(Z)(f) \frac{\psi(Z)(h_1)}{\lambda_0 h_0(x)} d\L_{x,e_0}(Z).\]
Let $\left(L_{x,\be}^{u,*}\right)_{x\in \XX,\be\in \E^{\ZZ_+}, u\in \UU}$ be an array of independent variables such that $L_{x,\be}^{u,*}\sim \L_{x,\be}^*$. We note once again $N_{x,\be}^{u,*}=\psi\left(L_{x,\be}^{u,*}\right)$. 
\\ If $u$ is at height $n$ in $\UU$, conditionally on $u\in \TT$, we note for short $L_u=L_{X(u),\xi_n}^{u}$ if $u\neq s_n$, and $L_{\s_n}=L^{\s_n,*}_{X(\s_n),\xi_n}$ the tuple valued variables describing the offspring of the individuals. Similarly, we note $N_u=\psi(L_u)$ their measure valued counterparts.
Then, conditionally on the environmental process $(\xi_n)$, we build recursively the biased random tree $(\TT^*,X^*)$, and the spine $(\s_n)_{n\geq 0}$ as follows.

We first set $X^*(\emptyset)=x$ and $\s_0=\emptyset$.  We define
\[\GG^*_{n+1}=\left\{ ui, u\in \GG^*_n\setminus \{\s_n\}, 1\leq i\leq N_u(\mathds{1}) \right\} \cup \left\{\s_n i,1\leq i\leq N_{\s_n}(\mathds{1})\right\}.\]
If $u\in \GG^*_n, u\neq \s_n$, conditionally on $N_{X(u),\xi_n}^u(\mathds{1})=d$ we set \[(X^*(u1), \cdots , X^*(ud))=L_u=L_{X(u),\xi_n}^u.\]
Moreover, conditionally on $N_{\s_n}(\mathds{1})=d$, \[(X^*(\s_n 1), \cdots , X^*(\s_n d))=L_{\s_n}=L^{\s_n,*}_{X(\s_n),\theta^n(\bxi)}.\] The next individual $\s_{n+1}$ of the spine is chosen at random among the children $\{\s_n 1,\dots \s_n d \}$ of $\s_n$ in such a way that $\s_{n+1}=\s_n i$ with probability proportional to $h_{n+1}(X^*(\s_n i)).$
Similarly as in the regular construction, we note $Z_n^*$ the occupation measure of the $n$-th generation of the labelled tree $(\TT^*,X^*)$. 
\[Z_n^*=\sum_{u\in \GG^*_n} \delta_{X^*(u)}.\]

In this subsection $\PP,\PP_{\bxi}$ refer to the respective annealed and quenched law of this whole construction.
Let us note $\mathcal{T}_n$ the $\sigma$-algebra on the set of labelled trees $\UU\times \XX^\UU$ generated by the events that only depend on the first $n$ generations of a tree, and $\mathcal{T}=\underset{n\geq 0} \bigvee \mathcal{T}_n$.
The labelled trees $(\TT,X)$ and $(\TT^*,X^*)$ are random variables on $(\UU\times \UU^\XX, \mathcal{T})$, we note $\P_{\bxi}$, $\P^*_{\bxi}$, $\P$ and $\P^*$ there respective quenched and annealed distributions. These are therefore probability measures on the space $(\UU\times \UU^\XX, \mathcal{T})$.
Note that $W_n$ is $\sigma(\bxi, \mathcal{T}_n)$-measurable. Moreover, the following analog of equation (5.3) of \cite{grama_rate_2014} holds
\begin{prop}
	For any $n\geq 0$, $\PP(d\bxi)$-almost surely,
	\[ d\P_{\bxi}^*|_{\mathcal{T}_n}=W_n d\P_{\bxi} |_{\mathcal{T}_n}.\]
\end{prop}
\begin{proof}
	For two non-labelled trees $t$ and $t\subset \UU$, we note $t\overset{n}{=}t'$ if the two trees coincide up to height $n$ and call $d_u(t)$ the number of offspring of the node $u\in t$.
	Let $t$ be a tree of height $n$, $(B_u)_{u\in t} \in \X^t$ be an array of measurable subsets of $\XX$ and $s\in \GG_n(t)$. We note $s=(i_1\dots i_n)$. Then the following equality of events holds :
	\begin{align}
		\left\{\TT^* \overset{n}{=} t \right\} \cap \bigcap_{u\in t} \left\{ X(u) \in B_u\right\} \cap \left\{ \s_n=s\right\}= \bigcap_{u \in t} \left\{ L_u \in \prod_{i=1}^{d_u(t)} B_u \right\} \cap \bigcap_{k=1}^n \left\{s_k=(i_1 \dots i_k)\right\} \nonumber
	\end{align}
	On this event, if $l(u)=d$, then conditionally on $\mathcal{T}_d$, either $L_u \sim \L_{X(u),\xi_d}$ or $L_{u}\sim \L_{X(u), T^d(\bxi)}^*$ depending on whether $u=\s_d$ or not. 
	This yields
	\begin{equation}
		\begin{split}&\PP\left[\left. \left\{\TT^* \overset{n}{=} t \right\} \cap \bigcap_{u\in t} \left\{ X(u) \in B_u\right\} \cap \left\{ \s_n=s\right\} \right| \bxi\right] \\&= \int \prod_{k=0}^{n-1} \prod_{\substack{u\in \GG_{k}(t)\\ u\neq (i_1\dots i_{k})}} \mathds{1}_{L_u\in \prod_{i=1}^{d_t(u)}B_{ui}} d\L_{X(u),\xi_{k}}(L_u) \\
			&\times \mathds{1}_{L_{k}\in \prod_{i=1}^{d_t(i_1\dots i_k)}B_{(i_1\dots i_k i)}} \frac{h_n(X(i_1\dots i_{k+1}))}{\sum_{i=1}^{d_{i_1\dots i_{k}}(t)} h_{k+1}(X(i_1\dots i_{k}i))} d\L^*_{X(i_1\dots i_k),T^{k}(\bxi)}(L_{(i_1\dots i_k)}) \end{split} \label{eq: loi arbre bias 1}
	\end{equation}
	In this context, the definition \eqref{eq:biased law} of the biased law $\L^*_{x,\be}$ yields, for all $k$,
	\begin{equation} \label{eq: csq def bias} d\L^*_{X(i_1\dots i_k),T^{k}(\bxi)}(L_{(i_1\dots i_k)}) =\frac{{\sum_{i=1}^{d_{i_1\dots i_{k}}(t)} h_{k+1}(X(i_1\dots i_{k}i))}}{\lambda_k h_k(X(i_1\dots i_k))}d\L_{X(i_1\dots i_k),\xi_k}(L_{(i_1\dots i_k)}).
	\end{equation}
	Therefore, plugging \eqref{eq: csq def bias} into \eqref{eq: loi arbre bias 1} yields
	\begin{align*}
		&\PP\left[\left. \left\{\TT^* \overset{n}{=} t \right\} \cap \bigcap_{u\in t} \left\{ X(u) \in B_u\right\} \cap \left\{ \s_n=s\right\} \right| \bxi\right] \\=& \int \prod_{k=0}^{n-1} \prod_{u\in \GG_{k}(t)} \mathds{1}_{L_u\in \prod_{i=1}^{d_t(u)}B_{ui}} d\L_{X(u),\xi_{k}}(L_u) 
		\times \frac{h_{n}(X(s))}{\lambda_{0,n} h_0(x)},
	\end{align*}
	where we recall that $x$ refers to the type of the initial individual $\emptyset=s_0$.
	Summing over all vertices $s$ in the $n$-th generation of $t$, we finally get
	\begin{align}
		\int \mathds{1}_{\left\{\TT \overset{n}{=} t \right\} \cap \bigcap_{u\in t} {X(u) \in B_u}}d\P^*_\bxi(\TT^*)&=\PP\left[\left. \left\{\TT^* \overset{n}{=} t \right\} \cap \bigcap_{u\in t} \left\{ X(u) \in B_u\right\} \right| \bxi\right] \nonumber \\ 
		&= \int \prod_{\substack{u\in t\\  l(u)\leq n-1}} \mathds{1}_{L_u\in \prod_{i=1}^{d_t(u)}B_{ui}} d\L_{X(u),\xi_{k}}(L_u) \times \frac{Z_n(h_n)}{\lambda_{0,n} h_0(x)} \nonumber \\
		&= \int \prod_{\substack{u\in t \\ l(u)\leq n-1}} \mathds{1}_{L_u\in \prod_{i=1}^{d_t(u)}B_{ui}} d\L_{X(u),\xi_{k}}(L_u) W_n \nonumber \\
		&= \int W_n \mathds{1}_{\left\{\TT \overset{n}{=} t \right\} \cap \bigcap_{u\in t} {X(u) \in B_u}} d\P_\bxi(\TT) \nonumber
	\end{align}
	The $\sigma$-algebra $\mathcal{T}_n$ is generated by events of the form $\left\{\TT \overset{n}{=} t \right\} \cap \bigcap_{u\in t} \{X(u) \in B_u\}$. Therefore the proposition is proved.
\end{proof}

\subsubsection{Nondegeneracy of the martingale}
We define in this section $W=\limsup W_n \in [0,\infty]$. Since $\PP(d\bxi)$-almost surely, the process $(W_n)$ is a $\P_\bxi$-nonnegative martingale, then $W=\lim W_n \in [0,\infty)$, $\P_\bxi$-almost surely, in $\PP(d\bxi)$-almost surely any environment sequence $\bxi$. Our main goal is to exhibit mild assumptions under which $W$ is not degenerate, in the sense that $\P_{\bxi}(W>0)>0, \PP(d\bxi)$-almost surely.
Let us fix a sequence $\be=(e_n)$. By Theorem 5.3.3 of \cite{durrett_probability_2010}, if $\P^*_\be(W<\infty)=1$ then $\int W d\P_\be=1$, and thus $\P_{\be}(W>0)>0$.
We focus now on the biased tree $\TT^*$ and prove that $\P^*_\bxi(\limsup W_n<\infty), \PP(d\bxi)-$a.s. 

To do so, following \cite{grama_kestenstigum_2023}, we introduce the quantity
\[A_n= \frac{(Z^*_n-\delta_{X^*(\s_n)})(h_n)}{\lambda_{0,n} Z^*_0(h_0)}-\sum_{k=0}^{n-1} \frac{(N_{\s_k}-\delta_{X^*(s_{k+1})})(h_{k+1})}{\lambda_{0,k+1}Z^*_0(h_0)}\]
We introduce $\Y=\sigma((\s_n,X^*(s_n),L_{\s_n})_{n\geq 0})$, the $\sigma-$algebra generated by the spine. Then the process $(A_n)$ satisfies
\begin{align} 
	\EE[A_{n+1}| \Y,\bxi,\mathcal{T}_n]&=\EE\left[\frac{(Z^*_{n+1}-\delta_{X^*(\s_{n+1})})(h_{n+1})-(N_{\s_n}-\delta_{X^*(\s_{n+1})})}{\lambda_{0,n+1} Z^*_0(h_0)} \right.\nonumber\\ & \left. \left. \qquad -  \sum_{k=0}^{n-1} \frac{(N_{\s_k}-\delta_{X^*(s_{k+1})})(h_{k+1})}{\lambda_{0,k+1}Z^*_0(h_0)} \right|\Y, \bxi, \mathcal{T}_n \right] \nonumber\\
	&=\frac{\EE[(Z^*_{n+1}-N_{\s_n})(h_{n+1}) | \bxi, \Y, \mathcal{T}_n]}{\lambda_{0,{n+1}}Z^*_0(h_0)}-\sum_{k=0}^{n-1} \frac{(N_{\s_k}-\delta_{X^*(s_{k+1})})(h_{k+1})}{\lambda_{0,k+1}Z^*_0(h_0)}.\label{eq:martingale A_1}\end{align}
Note that
\[Z^*_{n+1}= \sum_{u\in \GG^*_{n}\setminus \{\s_n\}} L_{\xi_n,X(u)}^u + N_{\s_n}.\]
Thus 
\begin{align}
	\EE[(Z^*_{n+1}-N_{\s_n})(h_{n+1}) | \bxi, \Y, \mathcal{T}_n]& =\sum_{u\in \GG^*_{n}\setminus \{\s_n\}} \delta_{X(u)}M_n h_{n+1} \nonumber \\
	&= (Z^*_n-\delta_{X^*(\s_n)}) M_n h_{n+1}\nonumber\\
	& = \lambda_n (Z^*_n-\delta_{X^*(\s_n)}) h_n \label{eq:martingale A_2}
\end{align} 
Plugging \eqref{eq:martingale A_2} into \eqref{eq:martingale A_1} proves that $(A_n)$ is a $\PP_{\bxi}$-martingale, conditionally on $\sigma(\bxi, \Y)$. 
To ensure the convergence of this martingale, we need the following statement.
\begin{lem}\label{lem:comportement geometrique immigration}
	We assume that there exists a nondecreasing and positive function $f$ such that $\frac{1}{x\log(x) f(x)}$ is integrable at infinity and 
	\[\EE \left[ \sup_{x\in \XX} \frac{ \EE\left[ Z_1(h_1) \log[ Z_1(h_1)] f(Z_1(h_1))| \bxi, Z_0=\delta_x \right] }{M_0h_1(x)}\right] <\infty.\]
	Then, $\PP$-almost surely, \[ \underset{n\rightarrow \infty}{\limsup}\frac{1}{n}\log^+ (N_{\s_n}(h_{n+1}))=0.\]
\end{lem}
Combining \autoref{ass:supercritical} with Proposition \ref{prop: prod lambda} and Lemma \ref{lem:comportement geometrique immigration} proves that the series 
\[ \sum_{k\geq 0} \frac{N_{\s_k}(h_{k+1})}{\lambda_{0,k+1} Z^*_0(h_0)}\]
is $\PP$-a.s. convergent. 
Moreover, noticing that $\delta_{X^*(s_{k+1}}(h_{k+1})\leq 1$, the series
\[\sum_{k\geq 0} \frac{(N_{\s_k}-\delta_{X^*(s_{k+1})})(h_{k+1})}{\lambda_{0,k+1}Z^*_0(h_0)}\]
also converges $\PP$-a.s. Its sum, which we denote $S$, is $\sigma(\bxi,\Y)$ measurable, thus the process $(A_n+S)_{n\geq 0}$ is a nonnegative martingale with respect to the filtration $\sigma(\Y,\bxi, \mathcal{T}_n)$. Thus, the sequence $(A_n)$ converges $\PP$-a.s. This proves in turn the $\PP$-almost sure convergence of 
\[\frac{(Z^*_n-\delta_{X^*(\s_n)})(h_n)}{\lambda_{0,n} Z^*_0(h_0)}. \]
Since, once again, the quantity $\delta_{X^*(\s_n)}(h_n)$ is bounded by $1$,
\[\frac{(Z^*_n)(h_n)}{\lambda_{0,n} Z^*_0(h_0)}\]
also converges $\PP$-a.s. to a finite variable.
As a consequence, for $\PP$-almost any sequence $\bxi=(\xi_n)$, the sequence $(W_n)_{n\geq 0}$ converges $\P^*_\bxi$-almost surely to a finite limit and $W=\limsup W_n<\infty$, $\P^*_\bxi$-almost surely. This concludes the proof.
\subsubsection{Proof of Lemma \ref{lem:comportement geometrique immigration}}
\begin{proof}
	Let us control first the tail of the annealed distribution of $N_{\s_n}(h_{n+1})$ :
	\begin{align}
		\PP[N_{\s_n}(h_{n+1})>t]&=\EE[\PP[N_{\s_n}(h_{n+1})>t|\bxi, X^*(\s_n)]] \nonumber\\
		&=\EE \left[\int \mathds{1}_{[t,+\infty)}(Z(h_{n+1})) d\L^*_{X^*(\s_n),\theta^n(\bxi)}(Z)\right] \nonumber\\
		&\leq \EE\left[\sup_{x\in\XX} \int \mathds{1}_{[t,+\infty)}(Z(h_{n+1})) d\L^*_{x,\theta^n(\bxi)}(Z)\right] \nonumber\\ 
		&\leq \EE\left[\sup_{x\in\XX} \int \mathds{1}_{[t,+\infty)}(Z(h_{1}))d\L^*_{x,\bxi}(Z)\right] \nonumber \\
		&\leq \EE\left[ \sup_{x\in\XX}\PP[ Z^*_1(h_1) >t | Z^*_0=x, \bxi]\right].\label{eq:controle tail 1}  
	\end{align}
	By the Markov inequality, we derive for any nonnegative and nondecreasing function $f$,
	\begin{equation} \label{eq:controle tail markov} \PP[ Z^*_1(h_1) >t | Z^*_0=x, \bxi]\leq \frac{\sup_{x\in\XX}\EE\left[ \log^+(Z^*_1(h_1)) f(Z_1^*(h_1)) | Z^*_0=\delta_x, \bxi=(e_n)\right]}{f(t) \log^+(t)}.\end{equation}
	Plugging \eqref{eq:controle tail markov} in \eqref{eq:controle tail 1}, we get
	\begin{align}
		\PP[N_{\s_n}(h_{n+1})>t]&\leq \EE\left[ \frac{\sup_{x\in\XX}\EE\left[ \log^+(Z^*_1(h_1)) f(Z_1^*(h_1)) | Z^*_0=\delta_x, \bxi=(e_n)\right]}{f(t) \log^+(t)}\right]\nonumber \\
		&\leq \frac{\EE\left[\sup_{x\in\XX}\EE\left[ \log^+(Z^*_1(h_1)) f(Z_1^*(h_1)) | Z^*_0=\delta_x, \bxi=(e_n)\right]\right]}{f(t) \log^+(t)}\label{eq:controle tail 2}
	\end{align}
	For $\PP$-almost any sequence $\be=(e_n)_{n\geq 0}$, setting $r(t)=t\log^+t f(t)$, we may write
	\begin{align}
		\frac{\EE\left[ r(Z_1(h_1)) |Z_0=\delta_x, \bxi=\be\right]}{M_0 h_1(x)}=&\int \log^+(Z_1(h_1)) f(Z_1(h_1)) \frac{Z_1(h_1)}{\lambda_0 h_0(x)} d\L_{x,e_0}(Z_1)\nonumber\\
		=&\int \log^+(Z(h_1)) f(Z(h_1)) d\L^*_{x,\be}(Z)\nonumber\\
		=&\EE\left[ \log^+(Z^*_1(h_1)) f(Z_1^*(h_1)) | Z_0=\delta_x, \bxi=(e_n)\right].\label{eq:lien biais}
	\end{align}
	Assumption \autoref{ass: LlogL} can be rewritten
	\[ \EE\left[ \sup_{x\in\XX} \frac{ \EE\left[ r(Z_1(h_1)) | Z_0=\delta_x, \bxi \right] }{M_0 h_1(x) } \right]<\infty\]
	which implies, by \eqref{eq:lien biais},
	\[ A:= \EE\left[\sup_{x\in\XX}\EE\left[ \log^+(Z^*_1(h_1)) f(Z_1^*(h_1)) | Z^*_0=\delta_x, \bxi=(e_n)\right]\right]<\infty.\]
	Hence, by \eqref{eq:controle tail 2}
	\[ \PP[N_{\s_n}(h_{n+1})>t]\leq \frac{A}{f(t)\log^+(t)}.\]
	Then, for any $c>0$,
	\begin{align} \EE\left[ \sum_{n\geq 1} \mathds{1}_{\frac{1}{n} \log( N_{\s_n}(h_{n+1}))>c}\right]&=\EE\left[ \sum_{n\in{\ZZ_+}} \mathds{1}_{\frac{1}{c} \log( N_{\s_n}(h_{n+1}))>n}\right]\nonumber \\
		&= \sum_{n\geq 1} \PP[\log^+( N_{\s_n}(h_{n+1}))>nc] \nonumber \\
		&= \sum_{n\geq 1} \PP[ N_{\s_n}(h_{n+1})>\exp(nc)] \nonumber \\
		&\leq \sum_{n\geq 1} \frac{A}{nc f(\exp{nc})}. \label{eq:borelcantelli immigration}
	\end{align}
	A elementary change of variable shows that since $\frac{1}{t\log(t) f(t)}$ is integrable at $+\infty$, then so is $\frac{1}{u f(\exp(u))}.$ Thus the sum \eqref{eq:borelcantelli immigration} converges, which proves that 
	\(\limsup_{n\rightarrow \infty} \frac{\log^+(N_{\s_n}(h_{n+1}))}{n} =0.\)
\end{proof}

\subsection{Almost sure convergence of the type distribution}
\label{subs:proof_type distrib}
{\color{black} Before proving Theorem \ref{thme:type_distrib}, let us discuss its statement.
	As pointed out in the introduction, we are not able to guarantee that the error term  ${o}\left(\left(\max(\tilde{\eta}{e^\lambda},e^{\frac{\lambda}{2}})+\varepsilon\right)^n\right)$ appearing in \eqref{eq: asymptotic development Z_n(f)} is negligeable with respect to the quantity $Z_0 M_{0,n}(f)$ for every $f\in \mathcal{B}(\XX)$. We highlight two specific subcases, which are neither exhaustive or exclusive of each other.
	\\ \textbf{Case 1:} If $\liminf_{n\rightarrow \infty} \left(\frac{Z_0 M_{0,n}(f)}{ \Vert Z_0 M_{0,n} \Vert_{TV}}\right)^{\frac{1}{n}}  > \max\left(\tilde{\eta},e^{-\frac{\lambda}{2}}\right)$, then conditionally on $\{W>0\}$ it holds
	\begin{equation} Z_n(f) \underset{n\rightarrow \infty}{\sim} W Z_0 M_{0,n}(f)\quad \PP\text{-a.s.}\label{eq:equivalent a.s. 1}\end{equation}
	\\ \textbf{Case 2:} If $A:=\limsup_{n\rightarrow \infty}  \left(\frac{Z_0 M_{0,n}(f)}{ \Vert Z_0 M_{0,n} \Vert_{TV}}\right)^{\frac{1}{n}}<1$, then it holds
	\begin{equation*} \limsup_{n\rightarrow \infty} \left(Z_n(f)^\frac{1}{n}\right)\leq \max(\tilde{\eta}{e^\lambda},e^{\frac{\lambda}{2}},A{e^\lambda})\quad \PP\text{-a.s.}\end{equation*}
	Any function $f$ such that $c\leq f\leq C$ for some positive constants $c,C$ (in particular $f=\mathds{1}$) belongs to case $1$ and satisfies \eqref{eq:equivalent a.s. 1}. This yields the first claim of Corollary \ref{cor:size and types}. However for functions $f$ such that $\inf_{x\in\XX}f(x)=0$, it seems possible that the typical distribution of types at time $n$ is concentrated on types associated with small values of $f$, in which case we might have, with positive probability $$\liminf_{n\rightarrow \infty} \left(\frac{Z_0 M_{0,n}(f)}{ \Vert Z_0 M_{0,n} \Vert_{TV}}\right)^{\frac{1}{n}}\leq\max(\tilde{\eta},e^{-\frac{\lambda}{2}}),$$ implying that along some subsequence $(n_k)$, the main term in the asymptotic expansion \eqref{eq: asymptotic development Z_n(f)} is in fact negligeable with respect ot the error term. Case 2 provides an upper bound on $Z_n(f)$ when one has additionally a sufficiently precise upper bound on the quenched mean of the distribution of types $\frac{Z_0 M_{0,n}(f)}{ \Vert Z_0 M_{0,n} \Vert_{TV}}$. 
	
	The claim \textit{ii)} of Theorem \ref{thme:type_distrib} provides a minorization condition relating a function $f$ and the coupling measure, which guarantees that $f$ belongs to Case 1. 
	
	Let us prove now Corollary \ref{cor:size and types} and then Theorem \ref{thme:type_distrib}
	\begin{proof}[Proof of Corollary \ref{cor:size and types}]
		Plugging $f=\mathds{1}$ in \eqref{eq: asymptotic development Z_n(f)} and choosing $\varepsilon$ small enough such that $\allowbreak \max(\tilde{\eta}{e^\lambda}, e^{\frac{\lambda}{2}})+\varepsilon<{e^\lambda}$, we obtain
		\[\Vert Z_n \Vert_{TV} =Z_n(\mathds{1})=W_n \Vert Z_0 M_{0,n}\Vert_{TV} +\underset{n\rightarrow \infty} {o}\left(e^{n\lambda}\right) \quad \PP\text{-a.s}.\]
		By Theorem \ref{thme: approx_KS}, \textit{iii)}, conditionally on $\{W>0\}$, we have 
		\[\lim_{n\rightarrow + \infty}\left(W_n \Vert Z_0 M_{0,n}\Vert_{TV}\right)^{\frac{1}{n}}={e^\lambda}\]
		and thus, conditionally on $\{W>0\}$,
		\[Z_n(\mathds{1})\sim W \Vert Z_0 M_{0,n}\Vert_{TV}  \quad \PP\text{-a.s}.\]
		Using again Theorem \ref{thme: approx_KS}, \textit{iii)}, we derive $\lim_{n\rightarrow \infty} \left(Z_n(\mathds{1})\right)^{\frac{1}{n}}={e^\lambda}.$ This proves \eqref{eq: size pop}.
		\\Let us choose now $f\in \mathcal{B}(\XX)$. We apply \eqref{eq: asymptotic development Z_n(f)} and derive, $\PP$-almost surely, conditionally on $\{W>0\}$, for any $\varepsilon>0$,
		\[\frac{Z_n(f)}{Z_n(\mathds{1})}=\frac{Z_0 M_{0,n}(f)}{\Vert Z_0 M_{0,n}\Vert_{TV}} \times \frac{W_n \Vert Z_0 M_{0,n}\Vert_{TV}}{Z_n(\mathds{1})}+ \underset{n\rightarrow\infty}{o}\left(\frac{\left(\max(\tilde{\eta}{e^\lambda},e^{\frac{\lambda}{2}})+\varepsilon\right)^n}{ Z_n(\mathds{1})}\right).\]
		By \eqref{eq: size pop}, we know that, if $\varepsilon>0$ is small enough, it holds $\PP$-almost surely 
		\[\lim_{n\rightarrow \infty} \frac{\left(\max(\tilde{\eta}{e^\lambda},e^{\frac{\lambda}{2}})+\varepsilon\right)^n}{ Z_n(\mathds{1})}=0\]
		conditionally on $\{W>0\}$, thus 
		\[\lim_{n\rightarrow +\infty} \frac{Z_n(f)}{Z_n(\mathds{1})}-\frac{Z_0 M_{0,n}(f)}{\Vert Z_0 M_{0,n}\Vert_{TV}} \times \frac{W_n \Vert Z_0 M_{0,n}\Vert_{TV}}{Z_n(\mathds{1})}=0,\]
		where $\frac{Z_0 M_{0,n}(f)}{\Vert Z_0 M_{0,n}\Vert_{TV}}$ is bounded by $1$. Hence, by \eqref{eq: size pop}, conditionally on $\{W>0\}$, it holds ${W_n \Vert Z_0 M_{0,n}\Vert_{TV}}\sim {Z_n(\mathds{1})}$. From this, we finally derive
		\[\lim_{n\rightarrow +\infty} \frac{Z_n(f)}{Z_n(\mathds{1})}-\frac{Z_0 M_{0,n}(f)}{\Vert Z_0 M_{0,n}\Vert_{TV}}=0.\]
\end{proof}}

\begin{proof}[Proof of Theorem \ref{thme:type_distrib}.]
	We note here $\pi_n=\frac{Z_0 M_{0,n}}{\Vert Z_0 M_{0,n}\Vert_{TV}}\in \mathcal{M}_{1}(\XX)$ the renormalized quenched mean of the population at time $n\geq 0$.
	Fix $f\in \mathcal{B}(\XX)$. We set \[f_n=f-\frac{\pi_n(f) h_n}{\pi_n(h_n)}=f-\frac{Z_0 M_{0,n}(f)}{Z_0M_{0,n}(h_n)}h_n=f-Z_0 M_{0,n}(f) \frac{h_n}{\lambda_{0,n}Z_0(h_0)}.\]
	Note that $\pi_n(f_n)=0$ and
	\begin{equation} \label{eq: norm f_n}
		\Vert f_n \Vert_\infty \leq \Vert f \Vert_{\infty} \left(1+\frac{1}{\pi_n(h_n)}\right)
	\end{equation}
	We also define 
	\begin{equation}
		\Delta_n=\frac{ Z_n(f_n)}{\lambda_{0,n}}=\frac{Z_n(f)}{\lambda_{0,n}}-W_n \frac{\pi_n(f)}{\pi_n(h_n)}Z_0(h_0)=\frac{Z_n(f)-W_n Z_0 M_{0,n}(f)}{\lambda_{0,n}}.
	\end{equation}
	The main idea of the proof is to show that the random sequence $(\Delta_n)_{n\geq 0}$ converges to $0$ exponentially fast in $\mathcal{L}^2(\PP_\bxi)$, $\as$ This will yield the desired result. 
	Notice first that
	\begin{align}
		Z_n(f_n)^2&=\left(\sum_{u\in\GG_n}f_n(X(u))\right)^2 \nonumber\\
		&=\sum_{u,v\in\GG_n}f_n(X(u))f_n(X(v)) \nonumber\\
		&=\sum_{u\in\GG_n}f_n(X(u))^2+\sum_{\substack{u,v\in\GG_n \\ u\neq v}}f_n(X(u))f_n(X(v)) \label{eq:somme paires}
	\end{align}
	where
	\begin{equation} \label{eq:somme paires egales}\sum_{u\in\GG_n}f_n(X(u))^2=Z_n(f_n^2).\end{equation}
	Moreover, each pair of distinct individuals $u\neq v \in\GG_n$ has a unique last common ancestor $w$ at some time $0\leq p \leq n-1$, characterized by the property that $w$ has two distinct children $u'\neq v'\in\GG_{p+1}$ such that $u'\preceq u$ and $v'\preceq v$. This can be used to reorganize the right hand side term in \eqref{eq:somme paires} as
	\begin{equation}
		\sum_{\substack{u,v\in\GG_n \\ u\neq v}}f_n(X(u))f_n(X(v))= \sum_{p=0}^{n-1}\underbrace{\sum_{w\in\GG_p} \sum_{\substack{u',v'\in\GG_{p+1}\\w\preceq u',v' \\ u'\neq v'}}\sum_{\substack{u\in\GG_n \\ u'\preceq u}} f_n(X(u))\sum_{\substack{v\in\GG_n \\ v'\preceq u}} f_n(X(v)).}_{:=\Sigma_p(f_n)}
		\label{eq:fourches} \end{equation}
	Combining \eqref{eq:somme paires}, \eqref{eq:somme paires egales} and \eqref{eq:fourches} yields
	\begin{equation}\label{eq:fourches complete}
		Z_n(f_n)^2=Z_n(f_n^2)+\sum_{p=0}^{n-1}\Sigma_p(f_n)
	\end{equation}
	We want now to compute and estimate the quenched expectation of \eqref{eq:fourches complete}. The left hand side term is the easiest to deal with :
	\begin{equation}\label{eq:estim first term}
		\EE_{\bxi}\left[Z_n(f_n^2)\right]=Z_0M_{0,n}(f_n^2)\leq \Vert f_n \Vert_{\infty}^2 \Vert Z_0 M_{0,n} \Vert_{TV}.
	\end{equation}
	The right hand side term requires a more subtle handling.
	For each pair of distinct individuals $u',v'$ at time $p+1$, the subtrees descending respectively form $u'$, $v'$ are independent conditionally on $\bxi$ and $\F_{p+1}$. This yields  
	\begin{equation}
		\EE_\bxi\left[\left.\Sigma_p(f_n) \right| \F_{p+1}\right]
		= \sum_{w\in\GG_p}\sum_{\substack{u',v'\in\GG_{p+1}\\w\preceq u',v' \\ u'\neq v'}} M_{p+1,n}(f_n)(X(u')) M_{p+1,n}(f_n)(X(v'))
		\label{eq:int fourches p+1}.
	\end{equation}
	Therefore, the quantity $\EE_{\bxi}\left[\Sigma_p(f_n)|\mathcal{F}_p\right]$ can be rewritten
	\begin{align*}
		\EE_{\bxi}\left[\Sigma_p(f_n)|\mathcal{F}_p\right]&=\EE_\bxi\left[\left. \EE_{\bxi}\left[\Sigma_p(f_n)|\mathcal{F}_{p+1}\right] \right| \mathcal{F}_p\right] \\
		&= \sum_{w\in\GG_p} \EE_{\bxi}\left[\left. \sum_{\substack{u',v'\in\GG_{p+1}\\w\preceq u',v' \\ u'\neq v'}} M_{p+1,n}(f_n)(X(u')) M_{p+1,n}(f_n)(X(v'))\right|\mathcal{F}_p  \right] \\
		&= \sum_{w\in\GG_p} V_p(M_{p+1,n}(f_n))(X(w)) \\
		&= Z_p(V_p(M_{p+1,n}(f_n))),
	\end{align*}
	where the quantity $V_p(g)(x)$ is defined for any $g\in\mathcal{B}(\XX)$ and $x\in\XX$ as
	\[V_p(g)(x)=\EE_\bxi\left[\left.\sum_{\substack{u,v\in\GG_{p+1} \\u\neq v}}g(X(u))g(X(v))\right|Z_p=\delta_x\right].\]
	Consequently,
	\[\EE_\bxi[\sigma_p(f_n)]=Z_0 M_{0,p }(V_p(M_{p+1,n}(f_n)))\]
	and
	
	\begin{align}
		\EE_{\bxi} \left[\sum_{p=0}^{n-1} \Sigma_p(f_n) \right]
		&= \sum_{p=0}^{n-1}Z_0 M_{0,p}(V_p(M_{p+1,n}(f_n))). \label{eq:somme fourche}\end{align}
	Notice that 
	$\left| V_p(g)\right| \leq V_p(\left|g\right|)$ and $V_p(\alpha g))=\alpha^2V_p(g)$ for any $g\in\B(\XX)$ and any $\alpha\in\RR$.
	Therefore, noting $\vvvert V_p \vvvert = \sup_{x\in\XX}V_p(\mathds{1})(x)$, it holds
	\begin{align}\label{eq: controle terme somme fourche}
		\left| Z_0 M_{0,p}(V_p(M_{p+1,n}(f_n)))\right|\leq \left(\sup_{x\in\XX} M_{p+1,n}(f_n)(x)\right)^2 \Vert Z_0 M_{0,p} \Vert \vvvert V_p \vvvert.
	\end{align}
	Using \eqref{eq: PF2} and recalling that $\pi_n(f_n)=0$, we derive
	\begin{align}
		\left| M_{p+1,n}(f_n)(x)\right|&= \left| M_{p+1,n}(f_n)(x)-\frac{h_k(x)}{\pi_0 M_{0,k}h_k} \Vert \pi_0 M_{0,n} \Vert_{TV} \pi_n(f_n)\right| \nonumber\\
		&\leq \Vert f_n \Vert_\infty \left\Vert \delta_x M_{p+1,n}-\frac{h_k(x)}{\pi_0 M_{0,k}h_k} \Vert \pi_0 M_{0,n} \Vert_{TV}\pi_n\right\Vert_{TV} \nonumber\\
		&\leq  \frac{4}{\gamma_{n-1}}\Vert f_n \Vert_\infty \Vert \delta_x M_{p+1,n} \Vert_{TV}\prod_{i=p+1}^{n-1} (1-\gamma_i) . \label{eq: application ergo fourches}
	\end{align}
	Plugging \eqref{eq: application ergo fourches} into \eqref{eq: controle terme somme fourche} yields
	\begin{equation}
		\left| Z_0 M_{0,p}(V_p(M_{p+1,n}(f_n)))\right|\leq \left(\frac{4}{\gamma_{n-1}}\Vert f_n \Vert_\infty \vvvert M_{p+1,n} \vvvert \prod_{p+1}^{n-1} (1-\gamma_i) \right)^2\left \Vert Z_0 M_{0,p} \right\Vert \vvvert V_p \vvvert.
	\end{equation}
	with
	\[\vvvert M_{p+1,n}\vvvert \leq \frac{1}{\gamma_p}\frac{\Vert Z_0 M_{0,n}\Vert_{TV}}{\Vert  Z_0M_{0,p+1}\Vert_{TV}},\]
	by \eqref{eq:coef gamma lower bound mult norm}. Thus
	\begin{align} &\left| Z_0 M_{0,p}(V_p(M_{p+1,n}(f_n)))\right|\nonumber\\ &\leq\left(\frac{4\Vert f_n \Vert_\infty}{\gamma_{n-1}\gamma_p} \frac{\Vert Z_0 M_{0,n}\Vert_{TV}}{\Vert  Z_0M_{0,p+1}\Vert_{TV}} \prod_{p+1}^{n-1} (1-\gamma_i) \right)^2\left \Vert Z_0 M_{0,p} \right\Vert_{TV} \vvvert V_p \vvvert
		\nonumber
		\\&\leq \underbrace{\left(\frac{4\Vert f_n \Vert_{\infty} \Vert Z_0 M_{0,n}\Vert_{TV}\prod_{i=0}^{n-1}(1-\gamma_i)}{\gamma_{n-1}}\right)^{2}}_{\alpha_n}
		\times \underbrace{\frac{\Vert Z_0 M_{0,p}\Vert_{TV} \vvvert V_p \vvvert}{\gamma_p\Vert Z_0M_{0,p+1}\Vert_{TV}^2\prod_{i=0}^{p}(1-\gamma_i)^2}}_{\beta_p}\label{eq:controle terme somme fourche 2}
	\end{align}
	Putting together \eqref{eq:fourches complete}, \eqref{eq:estim first term} with \eqref{eq:somme fourche} and \eqref{eq:controle terme somme fourche 2} yields, $\PP(d\bxi)$-a.s.
	\begin{equation*}
		\EE_{\bxi} \left[ Z_n(f_n)^2\right]\leq \Vert f_n \Vert_\infty^2 \Vert Z_0 M_{0,n}\Vert_{TV} + \alpha_n\sum_{p=0}^{n-1}\beta_p
	\end{equation*}
	and therefore
	\begin{equation}\label{eq:controle delta2}
		\EE_{\bxi}\left[\Delta_n^2\right]\leq \Vert f_n \Vert_{\infty}^2 \frac{\Vert Z_0 M_{0,n} \Vert}{\lambda_{0,n}^2}+\frac{\alpha_n}{\lambda_{0,n}^2}\sum_{p=0}^{n-1}\beta_p.
	\end{equation}
	To estimate the asymptotic behavior of the random sequences $\left(\left\Vert f_n \right\Vert_\infty\right)_{n\geq 0}$, $\left(\alpha_n\right)_{n\geq 0}$ and $\left(\beta_p\right)_{p\geq 0}$, we state the following lemma and postpone its proof to the end of the section.
	\begin{lem}\label{lem:estimates direction}
		Under assumptions \autoref{ass: boundedness_KS} to \autoref{ass:L2}, it holds $\PP(d\bxi)$-a.s.
		\begin{enumerate}[i)]
			\item $\limsup_{n\rightarrow \infty} \Vert f_n\Vert_\infty^\frac{1}{n}\leq1 $
			\item $\limsup_{n\rightarrow \infty} \left(\alpha_n\right)^{\frac{1}{n}}\leq \left(\tilde{\eta}{e^\lambda}\right)^2$
			\item $\displaystyle\limsup_{p\rightarrow \infty} \left(\beta_p\right)^{\frac{1}{p}} \leq\frac{{e^{-\lambda}}}{\tilde{\eta}^2}.$
		\end{enumerate}
	\end{lem}
	On the one hand, combining Lemma \ref{lem:estimates direction}, \textit{i)} with Theorem \ref{thme: approx_KS}, \textit{iv)} and Theorem \ref{prop: prod lambda}, we derive
	\begin{equation*}
		\lim_{n\rightarrow \infty}  \left( \Vert f_n \Vert \frac{\Vert Z_0 M_{0,n} \Vert}{\lambda_{0,n}^2}\right)^\frac{1}{n}\leq {e^{-\lambda}}<1 \quad \PP(d\bxi)\text{-a.s.}
	\end{equation*}
	which implies
	\begin{equation}\label{eq:controle first term}\Vert f_n \Vert \frac{\Vert Z_0 M_{0,n} \Vert}{\lambda_{0,n}^2}=\underset{n\rightarrow \infty}{O} \left(\left({e^{-\lambda}}+\varepsilon\right)^n\right) \quad \PP(d\bxi)\text{-a.s.} \end{equation}
	where ${e^{-\lambda}}+\varepsilon<1$ for $\varepsilon>0$ small enough.
	On the other hand, Lemma \ref{lem:estimates direction},\textit{iii)} clearly implies that $\as$, for any $\varepsilon>0$ and for any $p$ large enough, \( \beta_p\leq \left(\frac{{e^{-\lambda}}}{\tilde{\eta}^2}+\varepsilon\right)^p\). This implies 
	\begin{equation}\label{eq: estim somme fourche}
		\limsup_{n\rightarrow \infty}\left(\sum_{p=0}^{n-1} \beta_p \right)^\frac{1}{n}\leq \max\left(1,\frac{{e^{-\lambda}}}{\tilde{\eta}^2}\right) \quad \PP(d\bxi)\text{-a.s}.
	\end{equation}
	Indeed, 
	\begin{itemize}[-]
		\item If $\displaystyle\frac{{e^{-\lambda}}}{\tilde{\eta}^2}< 1$ then for $\varepsilon>0$ small enough $\displaystyle\frac{{e^{-\lambda}}}{\tilde{\eta}^2}+\varepsilon<1$ and
		\(\sum_{p=0}^{n-1}\beta_p=\underset{n\rightarrow\infty}{O}(1)\) $\PP(d\bxi)$-almost surely.
		\item If $\displaystyle\frac{{e^{-\lambda}}}{\tilde{\eta}^2}\geq 1$ then  $\displaystyle\frac{{e^{-\lambda}}}{\tilde{\eta}^2}+\varepsilon>1$ for all $\varepsilon>0$, thus
		\[\sum_{p=0}^{n-1}\beta_p=\underset{n\rightarrow \infty}{O} \left(\sum_{p=0}^{n-1}\left(\frac{1}{\tilde{\eta}{e^{2\lambda}}} +\varepsilon \right)^p \right)=\underset{n\rightarrow \infty}{O} \left(\left(\frac{1}{\tilde{\eta}{e^{2\lambda}}}+\varepsilon \right)^n \right) \quad \PP(d\bxi)\text{-a.s.}\]
		and taking an infimum over $\varepsilon>0$ yields \eqref{eq: estim somme fourche}. 
	\end{itemize} 
	Using \eqref{eq: estim somme fourche}, Lemma \ref{lem:estimates direction}, \textit{ii)} and Proposition \ref{prop: prod lambda}, we derive
	\begin{equation}
		\limsup_{n\rightarrow \infty}\left(\frac{\alpha_n}{\lambda_{0,n}^2}\sum_{p=0}^{n-1} \beta_p \right)^\frac{1}{n}\leq \max\left({e^{-\lambda}},\tilde{\eta}^2\right)<1 \quad \PP(d\bxi)\text{-a.s.}
	\end{equation}
	Thus for any $\varepsilon>0$,
	\begin{equation}\label{eq:controle final somme}
		\frac{\alpha_n}{\lambda_{0,n}^2}\sum_{p=0}^{n-1} \beta_p =\underset{n\rightarrow\infty}{O}\left(\left(\max\left({e^{-\lambda}},\tilde{\eta}^2\right)+\varepsilon\right)^n\right) \quad \PP(d\bxi)\text{-a.s.}
	\end{equation}
	Finally, plugging  \eqref{eq:controle first term} and \eqref{eq:controle final somme} into \eqref{eq:controle delta2}, we obtain that for any $\varepsilon>0$, it holds
	\begin{equation}
		\EE_\bxi[\Delta_n^2]=\underset{n\rightarrow\infty}{O}\left(\left(\max\left({e^{-\lambda}},\tilde{\eta}^2\right)+\varepsilon\right)^n\right)\underset{n\rightarrow\infty}{\longrightarrow}0 \quad \PP(d\bxi)\text{-a.s.}
	\end{equation}
	\color{black}
	Hence, if $\delta^2>\max\left({e^{-\lambda}},\tilde{\eta}^2\right)$, for any $\varepsilon>0$, it holds $\PP(d\bxi)-$almost surely
	\[\PP_\bxi[|\Delta_n|>\delta^n]\leq \frac{\EE_\bxi[\Delta^2_n]}{\delta^{2n}}=\underset{n\rightarrow \infty}{O}\left(\left(\frac{\max\left({e^{-\lambda}},\tilde{\eta}^2\right)+\varepsilon}{\delta^2}\right)^n\right)\]
	Consequently, by Borel-Cantelli's lemma, it holds $|\Delta_n|\leq \delta^n$ for $n$ large enough,  $\PP$-a.s.
	This proves that for all $\delta>\max\left({e^{-\frac{\lambda}{2}}},\tilde{\eta}\right)$,
	\[\frac{Z_n(f)}{\lambda_{0,n}}=W_n \frac{Z_0 M_{0,n}f}{\lambda_{0,n}}+\underset{n\rightarrow\infty}{O}(\delta^n) \quad \PP\text{-a.s.}\]
	Noticing that $\lim_{n\rightarrow \infty}\lambda_{0,n}^\frac{1}{n}={e^\lambda}$, we get finally, for any $\delta>\max(e^{\frac{\lambda}{2}}, \tilde{\eta}{e^\lambda})$, $\PP$-a.s., 
	\[Z_n(f)=W_n Z_0 M_{0,n} f+\underset{n\rightarrow +\infty}O(\delta^n).\]
	\color{black}
	Thus claim \textit{i)} is proved. Let us prove now \textit{ii)}.
	We assume now that the function $f\neq $ is nonnegative and satisfies $\EE\left[\log^{-}\left(\nu M_{1,N}(f)\right)\right]<+\infty$ for some deterministic integer $N\geq 1$.
	Let us write 
	\[\mu M_{0,n} f \geq \mu M_{0,n-N} M_{n-N,n-N+1}M_{n-N+1,n} f\]
	and notice that
	\[M_{n-N,n-N+1}M_{n-N+1,n}f(x)\geq c_{n-N} m_{n-N,n-N+1}(x) \nu_{n-N}(M_{n-N+1,n}f).\]
	Integrating the last inequality with respect to the measure $\mu M_{0,n-N}(dx)$, we obtain
	\begin{equation} \label{eq: cond f 0}\mu M_{0,n}f\geq c_{n-N}\Vert \mu M_{0,n-N+1}\Vert_{TV} \nu_{n-N}(M_{n-N+1,n}(f)). \end{equation}
	
	As a consequence from \cite[Lemma 3.3]{ligonniere_ergodic_2023}, we know that $\lim_{n\rightarrow \infty} \gamma_{n-N}^\frac{1}{n}=1$. Since $\gamma_{n-N}\leq c_{n-N} \leq 1$, we deduce \begin{equation} \lim_{n\rightarrow \infty }c_{n-N}^\frac{1}{n}=1. \label{eq: cond f 1} \end{equation}
	Finally, from the integrability property $\EE\left[\log^{-}\left(\nu M_{1,n} f\right)\right]<+\infty$, and the stationarity of the sequence $M_{n}$, we deduce that for all $\alpha>0$,
	$$ \sum_{n\geq 0} \PP\left[ \log^{-}\left(\nu_{n-N} M_{n-N+1,n}\right)\geq \alpha n \right] =\sum_{n\geq 0}\PP\left[ \log^{-}\left(\nu M_{1,N}\right)\geq \alpha n \right]<+\infty$$
	which yields 
	$$\sum_{n\geq 0} \PP\left[\left(\nu_{n-N} M_{n-N+1,n}\right)^\frac{1}{n} \leq e^{-\alpha}\right]<+\infty.$$
	Using Borel-Cantelli's lemma and letting $\alpha$ go to $0$, we deduce that 
	\begin{equation} \label{eq: cond f 2} \liminf_{n\geq 0} \left(\nu_{n-N} M_{n-N+1,n}\right)^\frac{1}{n} \geq 1 .\end{equation}
	Plugging \eqref{eq: cond f 1}, \eqref{eq: cond f 2} and \eqref{eq : growth rate_KS} into \eqref{eq: cond f 0} yields
	\[\liminf_{n\rightarrow \infty } \left( \mu M_{0,n}f\right)^\frac{1}{n}\geq {e^\lambda} \quad \as \]
	which can be strengthen into $\lim_{n\rightarrow \infty}  \left( \mu M_{0,n}f\right)^\frac{1}{n}={e^\lambda}$, using once again \eqref{eq : growth rate_KS} and the fact that $f$ is bounded.
	Therefore, we have  
	$\delta^n = o(Z_0 M_{0,n}(f))$ $\PP$-a.s, for any $\delta\in (\max(e^{\frac{\lambda}{2}}, \tilde{\eta}{e^\lambda}), {e^\lambda})$.
	Plugging this into \textit{i)} clearly yields both claims of \textit{ii)}.
	
	\color{black}
\end{proof}
Let us prove now Lemma \ref{lem:estimates direction}.
\begin{proof}[Proof of Lemma \ref{lem:estimates direction}.]
	First notice that, by \eqref{eq: norm f_n}, 
	\[\limsup_{n\rightarrow \infty}\left( \Vert f_n \Vert_{\infty}\right)^\frac{1}{n}\leq \limsup_{n\rightarrow \infty} \left( 1 +\frac{1}{\pi_n(h_n)} \right)^\frac{1}{n}.\]
	where, by definition of $\pi_n$ and by Lemma \ref{lem:eigenfunction}, it holds
	\[\pi_n(h_n)=\frac{Z_0 M_{0,n}h_n}{\Vert Z_0 M_{0,n}\Vert_{TV}}=Z_0(h_0)\frac{\lambda_{0,n}}{\Vert Z_0 M_{0,n} \Vert_{TV}}.\]
	Thus, applying Theorem \ref{thme: approx_KS}, \textit{iv)} and Lemma \ref{prop: prod lambda} we derive
	\(\lim_{n\rightarrow \infty}\left(\pi_n(h_n)\right)^\frac{1}{n}= 1,\)
	which in turns implies
	\(\limsup_{n\rightarrow \infty}\left( \Vert f_n \Vert_{\infty}\right)^\frac{1}{n}\leq 1, \as\)
	\textit{i)} is therefore proved.
	\\ By \cite[Lemma 3.3]{ligonniere_ergodic_2023}, it holds
	\begin{equation} \lim_{n\rightarrow \infty} \left(\prod_{i=0}^{n-1}(1-\gamma_i)\right)^\frac{1}{n}=\tilde{\eta} \text{ and } \lim_{n\rightarrow \infty} \left(\gamma_{n-1} \right)^\frac{1}{n}=1 \quad \as 
		\label{eq: estim prod gamma} 
	\end{equation}
	The estimate \textit{ii)} is now derived combining \textit{i)} with \eqref{eq: estim prod gamma} and Theorem \ref{thme: approx_KS}, \textit{iv)}.
	\\ Finally, notice that 
	\begin{align*}\left| V_p(g)(x)\right|&\leq\left| \EE_\bxi\left[\left.\sum_{\substack{u,v\in\GG_{p+1} \\u\neq v}}g(X(u))g(X(v))\right|Z_p=\delta_x\right]\right| \\
		& \leq \Vert g \vert_{\infty}^2 \left| \EE_\bxi\left[\left.\sum_{\substack{u,v\in\GG_{p+1} \\u\neq v}}1\right|Z_p=\delta_x\right]\right| \\
		&\leq \Vert g \vert_{\infty}^2 \left| \EE_\bxi\left[\left.\sum_{\substack{u,v\in\GG_{p+1}}}1\right|Z_p=\delta_x\right]\right|\\
		&\leq \vert g \Vert_\infty^2\EE_{\bxi}\left[\left. \Vert Z_{p+1}\Vert_{TV}^2 \right| Z_p=\delta_x \right]
	\end{align*}
	Thus $\vvvert V_p \vvvert\leq \sup_{x\in\XX} \EE_{\bxi}\left[\left. \Vert Z_{p+1}\Vert_{TV}^2 \right| Z_p=\delta_x \right]$. Moreover, assumption \autoref{ass:L2} ensures that the random quantity $\vvvert V_p\vvvert$ is $\PP(d\bxi)$-a.s. finite, and that 
	$\EE\left[\left|\log \vvvert V_p \vvvert\right|\right]<\infty.$
	This $\log$-integrability, combined with the stationarity of $(V_p)$ implies that $\lim_{p\rightarrow \infty}\left(\vvvert V_p\vvvert \right)^\frac{1}{p}=1$, by arguments similar to those of \cite[Lemma 3.3]{ligonniere_ergodic_2023}.
	Combining this with \eqref{eq: estim prod gamma} and Theorem \ref{thme: approx_KS}, \textit{iii)} yields 
	
	\begin{align*}\limsup_{p\rightarrow \infty} \left(\beta_p\right)^\frac{1}{p}& =\limsup_{p\rightarrow \infty} \left({\frac{\Vert Z_0 M_{0,p}\Vert_{TV} \vvvert V_p \vvvert}{\gamma_p\Vert Z_0 M_{0,p+1}\Vert_{TV}^2\prod_{i=0}^{p}(1-\gamma_i)^2}}\right)^\frac{1}{p} \\
		& \leq  \frac{\lim_{p\rightarrow \infty} \left( \vert Z_0 M_{0,p}\Vert_{TV}^\frac{1}{p}\right) \lim_{p\rightarrow \infty} \left( \vvvert V_p \vvvert^\frac{1}{p} \right)}{\lim_{p\rightarrow \infty} \gamma_p^\frac{1}{p} \left( \lim_{p\rightarrow \infty} \Vert Z_0 M_{0,p+1} \Vert ^\frac{1}{p}\right)^2 \left(\lim_{p\rightarrow \infty} \left(\prod_{i=0}^p (1-\gamma_i) \right)^{\frac{1}{p} }\right)^2} \\
		& \leq \frac{e^{-\lambda}}{{\tilde{\eta}^2}}.
	\end{align*} This concludes the proof of \textit{iii)}.
\end{proof}
\subsection{Extinction and explosion}
\label{subs:extinction-explosion}
\begin{proof}[Proof of Proposition \ref{prop:extinction-explosion}.]
	Let us consider a MGWRE and assume that conditionally on survival, it holds almost surely $\limsup_{n\rightarrow \infty}Z_n(\{x_0\})=+\infty$. In other words,
	\begin{equation}\label{eq:survie type 2}\PP_{\bxi}[\ext]+\PP_{\bxi}[\limsup Z_n({x_0})=+\infty]=1\quad \PP[d\bxi]\text{-a.s}\end{equation}
	We additionally assume that assumptions \autoref{ass: boundedness_KS} to \autoref{ass: LlogL} hold, therefore $\PP_{\bxi}[W>0|Z_0=\delta_x]>0, \as$, for any type $x\in\XX$ and in particular for $x=x_0$.
	
	Let us introduce the stopping time
	\(T_K=\inf\{n\geq 1| Z_n(x_0)\geq K\}\)
	and note for each $(x,\be)\in\XX\times \E^{\ZZ_+}$,
	\[q_{\nexp}(x,\be)=\PP[\nexp|Z_0=\delta_x,\bxi=\bar{e}]=\PP[W>0|Z_0=\delta_x, \bxi=\be]\leq 1.\]
	Since 
	\(\PP_{\bxi}[\nexp|Z_n]=\prod_{x\in\XX} q_{\nexp}(x,\theta^n(\bxi))^{Z_n(x)},\)
	we have
	\[\PP_{\bxi}[\nexp|T_K=n]\leq q_{\nexp}({x_0},\theta^n(\bxi))^{K}.\] Thus, 
	\[\PP_{\bxi}[\nexp\cap \{T_K<\infty\}]\leq \sum_{n\geq 1} \PP_{\bxi}[T_K=n]q_{\nexp}({x_0},\theta^n(\bxi))^{K},\]
	where $\allowbreak\PP_{\bxi}[T_K=n]$ is $\sigma(\xi_0,\dots \xi_{n-1})$-measurable and $q_{\nexp}({x_0},\theta^n(\bxi))$ is $\sigma((\xi_k)_{k\geq n})$-measurable. Since the sequence $(\xi_n)$ is here assumed to be i.i.d, the quantities $\PP_{\bxi}[T_K=n]$ and $q_{\nexp}({x_0},\theta^n(\bxi))$ are independent from each other. Hence, integrating with respect to the environment, we get 
	\begin{align}
		\PP[\nexp\cap \{T_K<\infty\}] &\leq \sum_{n\geq 1} \PP[T_K=n]\EE \left[ q_{\nexp}({x_0},\theta^n(\bxi))^{K}\right] \nonumber \\
		&=\sum_{n\geq 1} \PP[T_K=n]\EE \left[ q_{\nexp}({x_0},\bxi)^{K}\right] \nonumber \\
		&=\PP[T_K<\infty]\EE \left[ q_{\nexp}({x_0},\bxi)^{K}\right].
	\end{align}
	By Theorem \ref{thme:nondegenerescence},  $q_{\nexp}({x_0},\bxi)<1 \, \as$ This yields
	\( \lim_{K\rightarrow \infty} q_{\nexp}({x_0},\bxi)^K=0,{\PP(d\bxi)}\text{-a.s.}\)
	and thus \(\lim_{K\rightarrow \infty} \PP[\nexp\cap \{T_K<\infty\}] =0.\)
	We derive therefore
	\begin{align*}\PP\left[\nexp \cap \{\limsup{Z_n(0)}=+\infty\}\right]&=\PP\left[\nexp \cap \bigcap_{K\geq 1}\{T_K=\infty\}\right]\\
		&=\lim_{K\rightarrow \infty} \PP[\nexp\cap \{T_K<\infty\}] \\& =0,
	\end{align*}
	hence
	\[\PP_{\bxi}\left[\nexp \cap \{ \limsup{Z_n(0)}=+\infty\}\right]=0 \quad \PP(d\bxi)\text{-a.s.}\]
	Finally, the equality
	\[\PP\left[\nexp \cap \{\limsup {Z_n(0)}=+\infty\}\right]=\EE\left[\PP_{\bxi}\left[\nexp \cap \{\limsup{Z_n(0)}= +\infty\}\right]\right]\]
	yields
	\[\PP_{\bxi}\left[\nexp \cap \{\limsup{Z_n(0)}= +\infty\}\right]=0 \quad \PP(d\bxi)\text{-a.s.}\]
	From this and \eqref{eq:survie type 2}, we derive $\as$,
	\[\PP_{\bxi}\left[\nexp \right] + \PP_{\bxi}\left[\{\limsup{Z_n(0)}= +\infty\}\right]\leq 1=\PP_{\bxi}\left[\ext \right] + \PP_{\bxi}\left[\{\limsup{Z_n(0)}= +\infty\}\right].\]
	Hence $\PP_{\bxi}\left[\nexp \right]\leq \PP_{\bxi}\left[\ext \right], \as$ Combining this with the clear inclusion $\ext \subset \nexp$ proves that $\ext$ and $\nexp$ only differ by a $\PP$-negligible event.
\end{proof}
\section{Discrete time model for an age structured population}
\label{sec:LGWRE}
\subsection{Presentation of the model}
In this subsection, we focus on an example of population model where the type accounts for the age of the individuals. We choose to model this notion of the age of an individual by the integer number of generations since its birth. Therefore, we take $\XX={\ZZ_+}$, and $\X=\mathcal{P}({\ZZ_+})$. At each time step, an individual of time $x$ may create a random number of individuals of age $0$, which we refer to as \textit{newborns}. Moreover, at each time step,  after giving birth to newborns, each individual may die or survive to the next time step, in which case its offspring will additionally contain an individual of age $x+1$.
For each $x\in \XX$, $e\in \E$, let $F_{x,e}$ be an integer valued variable of mean $f_{x,e}$ (encoding the number of newborns), and $S_{x,e}$ be a Bernoulli variable of parameter $s_{x,e}$ (encoding the survival or death of the individual), independent of $F_{x,e}$. We define $L_{x,e}=(0,\dots, 0)\in \XX^{k}$ if $(S_{x,e},L_{x,e})=(0,k)$ and $L_{x,e}=(x+1,0,\dots, 0)\in \XX^{k+1}$ if $(S_{x,e},L_{x,e})=(1,k)$. Choosing $\L_{x,e}$ as the distribution of $L_{x,e}$, we introduce a particular class of GWRE which is adapted to model an age structured population.
In particular, we write $N_{x,e}^{k,n}=\psi(L_{x,e}^{k,n})= F_{x,e}\delta_0 + S_{x,e}\delta_{x+1}.$
Noting $f_{x,e}=\EE[F_{x,e}]$ and $s_{x,e}=\EE[F_{x,e}]$, the mean matrix of such a process is of the form
\begin{equation*} M(e)=\begin{pmatrix} \label{eq : Leslie_KS}
		f_{0,} & s_{0,e} &  0  &  0  & \dots  \\
		f_{1,e} &  0  & s_{1,e} &  0  & \dots  \\
		f_{2,e} &  0  &  0  & s_{2,e} & \ddots \\
		f_{3,e} &  0  &  0  &  0  & \ddots \\
		\vdots & \vdots & \vdots & \vdots & \ddots \\
	\end{pmatrix}. \end{equation*}
This is clearly an infinite dimensional generalization of the Leslie matrices presented for example in \cite{caswell_life_2010}. For this reason we call a MGWRE with such offspring distribution a Leslie-GWRE or (L-GWRE for short). 
In this context, the random variables $Z_n$ are point measures on the discrete type space $\XX=\ZZ_+$. Therefore $Z_n$ is characterized by the vector $(Z_n(k))_{k\geq 0}$, where we note for short $Z_n(k)=Z_n(\{k\})$.
\subsection{Preliminary assumption}
We recall that if $X_1$ and $X_2$ are two real valued random variables we say that $X_2$ dominates stochastically $X_1$ (and note $X_1\preceq X_2$) if 
\[\forall t\in\RR,\, \PP[X_1\geq t] \leq \PP[X_2\geq t].\]
In the whole section, we assume that the following assumption holds.
\begin{ass}\label{ass:leslie} The random variables $(F_{x,e},S_{x,e})_{x\in\XX,e\in\E}$ satisfy
	\begin{enumerate}[i)]
		\item $s_{x,e},f_{x,e}>0$ for any $x\in\XX$, $e\in\EE$.
		\item  $F_{x+1,e}\preceq F_{x,e}$ and $S_{x+1,e}\preceq S_{x,e}$ for any $x\in \XX, e\in\EE$.
		\item $\EE[\log f_{0,\xi_0}]<\infty$
		\item $\EE\left[\left|\log\left(\sup_{x\in\XX} \frac{s_{x,\xi_0}}{f_{x,\xi_0}}\right)\right|\right]<\infty$
	\end{enumerate}
\end{ass}
By \cite[Prop 5.1]{ligonniere_ergodic_2023}, we have

\begin{prop}\label{prop:hyp_leslie_2}
	Consider a L-GWRE $(Z_n)_{n\geq 0}$, with a stationary and ergodic environmental sequence $\bxi=(\xi_n)_{n\geq 0}$, such that the offspring random variables $(F_{x,e},S_{x,e})_{x\in\XX,e\in \E}$ satisfy \autoref{ass:leslie}.
	Then the associated quenched first moment semi-group $(M_n)$ satisfies Assumptions \autoref{ass: boundedness_KS} to \autoref{ass : moments_cd_KS}, and therefore the conclusions of Theorems \ref{thme: approx_KS} and \ref{thme: martingale} hold.
\end{prop}
In \autoref{ass:leslie}, we restrict ourselves  to offspring distributions $(F_{x,e})_{x\geq 0},(S_{x,e})_{x\geq 0}$ which are nonincreasing as a function of $x$, with respect to the stochastic domination relation. This is consistent with the intuition that older individuals have both a higher chance of dying and tend to make less children than younger ones.
\begin{remark}
	In some biological contexts, one might expect that the peak age in terms of fertility or survival is some age $x_0>0$. This is not compatible with Assumption \autoref{ass:leslie}-\textit{ii)}. To allow this situation to happen, one can weaken \textit{ii)} and assume instead that the variables $(F_{x,e})$ and $(S_{x,e})$ are decreasing only after some rank $x=x_0$. Then some additional controls on the variables $F_{x,e}$ and $S_{x,e}$ for $x\leq x_0$ allow to deduce \autoref{ass: boundedness_KS}-\autoref{ass : moments_cd_KS} using \cite[Prop 5.1, Remark 5.2]{ligonniere_ergodic_2023}.
\end{remark}
\subsection{Non degeneneracy of $(W_n)$ for L-GWRE} \label{subs : ex non deg}
Let us now focus on assumption \autoref{ass: LlogL} and see how to check it when working with L-GWREs.

\begin{lem}\label{lem:LlogL for leslie}
	Let $(Z_n)$ be a L-GWRE satisfying \autoref{ass:leslie} and assume that there exists some $\varepsilon>0$ such that at least one of the following assertions holds :
	\begin{equation} \label{eq:LlogL age}  \EE\left[ \sup_{x\in \XX}  \frac{\sum_{k\geq 0}(k+1)\log^+(k+1)^{1+\varepsilon}\PP\left[F_{x,\xi_0}= k|\bxi\right]}{f_{x,\xi_0}}\right] <\infty, \end{equation}
	
	\begin{equation} \label{eq:LlogL age tail}  \EE\left[\sup_{x\in \XX}  \frac{\sum_{k\geq 0}\log^+(k+1)^{1+\varepsilon} \PP\left[F_{x,\xi_0}\geq k|\bxi \right]}{f_{x,\xi_0}}\right] <\infty. \end{equation}
	Then \autoref{ass: LlogL} holds. If additionally $\lambda>0$, then the conclusions of Theorem \ref{thme:nondegenerescence} hold.
\end{lem}
We postpone the proof of this lemma at the end of Subsection \ref{subs : ex non deg}. This statement proves that for L-GWREs, assumption \autoref{ass: LlogL} is a consequence of some uniform integrability assumption on the variables $(F_{x,e})_{x,e}$. Unfortunately, the argument of the supremum both in \eqref{eq:LlogL age} and \eqref{eq:LlogL age tail} is not monotonous in $F_{x,e}$, thus the monotonicity assumption \autoref{ass:leslie}-\textit{ii)} does not allow to simply remove the supremum in criterions \eqref{eq:LlogL age} and \eqref{eq:LlogL age tail}. However, these criterions can be pretty easily checked assuming some control on the tails of the variables $(F_{x,e})_{x,e}$. For any real valued sequences $(u_k)$ and $(v_k)$ and any $\beta>0$, we set $$ u_k\overset{\beta}{\asymp} v_k \Leftrightarrow \beta^{-1} v_k \leq u_k \leq \beta v_k \text{ for all } k\geq 0.$$ Let us use this notation to introduce three classes of distributions for $(F_{x,e})_{x\in\XX, e\in\E}$.
\begin{defi}
	Let $(F_{x,e})_{x\in\XX, e\in\E}$ be a familiy of integer-valued random variables. We say that $(F_{x,e})_{x\in\XX, e\in\E}$ belongs to
	\begin{itemize}[-]
		\item \textbf{BS} if and only if for each $(x,e)\in\XX\times \E$, there exists a deterministic integer $A_{x,e}$ and a real number $\varepsilon >0$ such that 
		\begin{align*} \PP[F_{x,e}\geq A_{x,e}]=0 \text{ and } \EE\left[\sup_{x\in\XX} \frac{A_{x,\xi_0}\left(\log^+(A_{x,\xi_0})\right)^{1+\varepsilon}}{f_{x,\xi_0}}\right]<\infty,
		\end{align*}
		\item \textbf{ExpTail} if and only if for each $(x,e)\in\XX\times \E$, there exists $q(x,e)<1, \alpha(x,e)>0$ and $\beta(e)>0$ such that
		$$\PP[F_{x,e}\geq k] \overset{\beta(e)}{\asymp}\alpha(x,e) q(x,e)^k$$ and
		\begin{align*}
			\EE\left[\frac{\beta(\xi_0)^2}{1-\sup_{x\in\XX}q(x,\xi_0)}\right]&<\infty.
		\end{align*}
		\item \textbf{PolyTail} if and only if for each $(x,e)\in\XX\times \E$, there exists $\delta(x,e)>1, \alpha(x,e)>0$ and $\beta(e)>0$, it holds
		$$\PP[F_{x,e}\geq k] \overset{\beta(e)}{\asymp}\alpha(x,e) k^{-\delta(x,e)},$$ where 
		\begin{align*}
			\EE\left[\beta(\xi_0)^2 \inf_{x\in \XX} \delta(x,\xi_0)\right]<\infty \text{ and } \EE\left[\frac{\beta(\xi_0)^2}{(\inf_{x\in \XX} \delta(x,\xi_0)-1)^{1+\varepsilon}}\right]&<\infty \text{ for some $\varepsilon>0$.}
		\end{align*}
	\end{itemize}
\end{defi}
When the random variables $(F_{x,e})_{x\in\XX, e\in\E}$ are in one of the three previous classes, our assumption \autoref{ass: LlogL} can be verified, allowing to apply Theorem \ref{thme:nondegenerescence}.
{\color{black}
	\begin{prop}
		Consider a supercritical L-GWRE satisfying \autoref{ass:leslie}. We assume that the family of variables $(F_{x,e})_{x,e}$ belongs to one of the classes \textbf{BS}, \textbf{ExpTail}, \textbf{PolyTail}. Then, the process $(W_n)$ is well defined, it is a nonnegative martingale and its limit $W$ satisfies\[\PP_\bxi[W>0]>0\quad\PP(d\bxi)\text{-a.s.}\]
\end{prop}}
\begin{proof}
	by Proposition \ref{prop:hyp_leslie_2}, assumptions \autoref{ass: boundedness_KS} to \autoref{ass : moments_cd_KS} hold when \autoref{ass:leslie} is satisfied. Since we additionally assume that $\lambda>0$, it only remains to check \autoref{ass: LlogL} to apply Theorem \ref{thme:nondegenerescence}. Let us check it in each of the three cases.
	\begin{itemize}[-]
		\item If $(F_{x,e})_{x\in\XX, e\in\E}$ belongs to \textbf{BS}, then
		{		\color{black}\[\sum_{k\geq 0} (k+1)\log(k+1)^{1+\varepsilon}  \frac{\PP\left[F_{x,e}=k\right]}{f_{x,e}}\leq \frac{(A_{x,e})\log(A_{x,e})^{1+\varepsilon}}{f_{x,e}}.\]}
		Therefore, applying criterion \eqref{eq:LlogL age}, we can show that \autoref{ass: LlogL} is satisfied as soon as 
		\[ \EE\left[\sup_{x\in\XX} \frac{A_{x,\xi_0}\log(A_{x,_{\xi_0}})^{1+\varepsilon}}{f_{x,\xi_0}}\right]<\infty.\]
		\item Assume now that $(F_{x,e})_{x\in\XX,e\in \E}$ belongs to \textbf{ExpTail}. In this case, we apply criterion \eqref{eq:LlogL age tail}.
		On the one hand, since $\PP[F_{x,e}\geq k]\geq \beta(e)^{-1} \alpha(x,e) q(x,e)^k$, we can derive
		\begin{align*}
			f_{x,e}&= \EE[F_{x,e}] \\
			&=\sum_{k\geq 1} \PP[F_{x,e}\geq k] \\
			&\geq \beta(e)^{-1} \alpha(x,e) \sum_{k\geq 1} q(x,e)^k\\
			&\geq \beta(e)^{-1} \alpha(x,e) \frac{q(x,e)}{1-q(x,e)}.\\
		\end{align*}
		On the other hand, for any $\varepsilon>0$,
		\begin{align*}
			\sum_{k\geq 0}\log^+(k+1)^{1+\varepsilon} \PP\left[F_{x,e}\geq k\right] &\leq \beta(e)\alpha(x,e)\sum_{k\geq 1} \log^+(k+1)^{1+\varepsilon} q(x,e)^k.\\ 
		\end{align*}
		Consequently, noting $C=\sup_{k\geq 1} \frac{\log^+(k+1)^{1+\varepsilon}}{k}$,
		\begin{align*}
			\frac{\sum_{k\geq 0}\log^+(k+1)^{1+\varepsilon} \PP\left[F_{x,e}\geq k\right]}{f_{x,e}}&\leq \beta(e)^2 (1-q(x,e))\sum_{k\geq 1}\log^+(k+1)^{1+\varepsilon} q(x,e)^{k-1}\\
			&\leq C\beta(e)^2(1-q(x,e)) \sum_{k\geq 1}kq(x,e)^{k-1} \\
			&\leq C\beta(e)^2(1-q(x,e)) \frac{1}{(1-q(x,e))^2} \\
			&\leq C\frac{\beta(e)^2}{1-q(x,e)}.
		\end{align*}
		Taking a supremum over $x\in \XX$ and integrating over the environmental space, we get indeed
		\begin{equation*}
			\EE\left[ \sup_{x\in\XX} \frac{\sum_{k\geq 0}\log^+(k+1)^{1+\varepsilon} \PP\left[F_{x,\xi_0}\geq k|\bxi\right]}{f_{x,e}}\right]\leq C(\varepsilon)\EE\left[\frac{\beta(\xi_0)^2}{1-\sup_{x\in\XX} q(x,\xi_0)}\right]<\infty,
		\end{equation*}
		which proves \autoref{ass: LlogL}.
		\item We suppose finally that $(F_{x,e})_{x\in\XX,e\in \E}$ belongs to \textbf{Polytail}.
		On the one hand, we derive similarly as in the previous case
		\begin{align*}
			f_{x,e}	&=\sum_{k\geq 1} \PP[F_{x,e}\geq k] \\
			&\geq \beta(e)^{-1} \alpha(x,e) \sum_{k\geq 1} k^{-\delta(x,e)}\\
			&\geq \beta(e)^{-1}\alpha(x,e) \int_1^{+\infty} t^{-\delta(x,e)}dt \\
			&\geq \frac{\alpha(x,e)}{(\delta(x,e)-1)\beta(e)}
		\end{align*}
		and for any $\varepsilon>0$,
		\begin{align*}
			\sum_{k\geq 0}\log^+(k+1)^{1+\varepsilon} \PP\left[F_{x,e}\geq k\right] &\leq \beta(e)\alpha(x,e)\sum_{k\geq 1} \log^+(k+1)^{1+\varepsilon} k^{-\delta(x,e)}. 
		\end{align*}
		Hence, we get
		\begin{align*}
			\frac{\sum_{k\geq 0}\log^+(k+1)^{1+\varepsilon} \PP\left[F_{x,e}\geq k\right]}{f_{x,e}}&\leq \beta(e)^2 (\delta(x,e)-1)\sum_{k\geq 1}\log^+(k+1)^{1+\varepsilon} k^{-\delta(x,e)}\\
			&\leq \beta(e)^2(\delta(x,e)-1) u(1+\varepsilon,\delta(x,e))\end{align*}
		where $u:(s,a)\in(1,\infty)^2\mapsto \sum_{k\geq 1} \log(k+1)^s k^{-a}.$
		A careful study of $u$ (see Appendix, Lemma \ref{lem: etude u}) allows to show that 
		\[u(1+\varepsilon,\delta(x,e))\leq C_1+\frac{C_2}{(\delta(x,e)-1)^{2+\varepsilon}},\]
		where $C_1$ and $C_2$ are two explicit numbers depending only on $\varepsilon$.
		This yields
		\begin{align*}
			\frac{\sum_{k\geq 0}\log^+(k+1)^{1+\varepsilon} \PP\left[F_{x,e}\geq k\right]}{f_{x,e}}&\leq C_1 \beta(e)^2(\delta(x,e)-1) \\&+
			C_2 \frac{\beta(e)^2}{(\delta(x,e)-1)^{1+\varepsilon}}
		\end{align*}
		and finally 
		\begin{align*}
			\EE\left[ \sup_{x\in\XX} \frac{\sum_{k\geq 0}\log^+(k+1)^{1+\varepsilon} \PP\left[F_{x,\xi_0}\geq k|\bxi\right]}{f_{x,e}}\right]&\leq  C_1 \EE\left[\beta(\xi_0)^2(\inf_{x\in\XX}\delta(x,\xi_0)-1)\right] \\ &+
			C_2 \EE\left[\frac{\beta(\xi_0)^2}{(\inf_{x\in\XX}\delta(x,\xi_0)-1)^{1+\varepsilon}}\right]<\infty.
		\end{align*}
		As a consequence, \autoref{ass: LlogL} also holds on the class \textbf{PolyTail}.
		
	\end{itemize}
\end{proof}
\begin{proof}[Proof of Lemma \ref{lem:LlogL for leslie}]
	We assume that the distributions of the variables $(F_{x,e},S_{x,e})$ satisfy \autoref{ass:leslie}.
	Let us check that \eqref{eq:LlogL age} and \eqref{eq:LlogL age tail} both imply \autoref{ass: LlogL}.
	Conditionally on $Z_0=\delta_x$, it holds 
	\[ h_1(0)F_{x,\xi_0}\leq Z_1(h_1)=S_{x,\xi_0}h_1(x+1)+F_{x,\xi_0} h_1(0)\leq (1+F_{x,\xi_0})\]
	since $S_{x,\xi_0}\in \{0,1\}$ and $h_1$ is positive and bounded by $1$.
	Therefore
	\begin{equation} \label{eq: moyenne taille_0} M_0h_1(x)=\EE[Z_1(h_1) | \bxi, Z_0=\delta_x ]\geq f_{x,\xi_0}h_1(0).\end{equation}
	Moreover, notice that $h_1(0)=\lim_{n\rightarrow \infty}\frac{ \Vert\delta_0 M_{1,n}\Vert_{TV}}{\vvvert M_{1,n} \vvvert} \geq d_1$
	by definition of $d_1$. By Lemma 5.5 and Remark 5.2 of \cite{ligonniere_ergodic_2023}, we can choose $d_1$ in such a way that 
	\[1\geq d_1\geq \inf_{x,n} \frac{s_{0,\xi_0}\dots s_{n,\xi_n}f_{n+1,\xi_{n+1}}}{s_{x+0,\xi_0} \dots s_{x+n,\xi_n}f_{x+n+1,\xi_{n+1}}}\]
	Moreover, \autoref{ass:leslie}-\textit{ii)} implies that the sequences $(f_{x,e})_{x\geq 0}$ and $(s_{x,e})_{x\geq 0}$ are nonincreasing, thus 
	\[\inf_{x,n} \frac{s_{0,\xi_0}\dots s_{n,\xi_n}f_{n+1,\xi_{n+1}}}{s_{x+0,\xi_0} \dots s_{x+n,\xi_n}f_{x+n+1,\xi_{n+1}}}\geq 1\]
	and $d_1=1$. Therefore, under \autoref{ass:leslie}, it holds $h_1(0)\geq 1$ $\PP(d\xi)$-a.s., hence $h_1(0)=0$ since $\Vert h_1 \Vert_{\infty}=1$.
	Plugging this in \eqref{eq: moyenne taille_0} yields
	\begin{equation} \label{eq: moyenne taille} M_0h_1(x)\geq f_{x,\xi_0}.\end{equation}
	We choose to verify \autoref{ass: LlogL} with $f(t)=\log^+(t)^{\varepsilon}$, which satisfies the integrability assumption $\int^{+\infty} (t\log^+(t) f(t))^{-1}dt<\infty$. Therefore, we set $g(t)=t(\log^+(t))^{1+\varepsilon}$ and derive
	\begin{align}
		\EE\left[ g(Z_1(h_1))| Z_0=\delta_x, \bxi \right]&\leq \EE\left[g(F_{x,\xi_0}+1)| \bxi\right] \nonumber \\
		&\leq \sum_{k\geq 0} g(k+1) \PP\left[F_{x,\xi_0}=k|\bxi \right] \nonumber \\
		& \leq \sum_{k\geq 0} (k+1) \log^+(k+1)^{1+\varepsilon} \PP\left[F_{x,\xi_0}=k|\bxi \right] \label{eq: LogL avant IPP}.
	\end{align}
	Putting together \eqref{eq: LogL avant IPP} and \eqref{eq: moyenne taille}, taking a supremum in $x$ and integrating with respect to $\bxi$, we show that
	\[\EE \left[ \sup_{x\in \XX} \frac{ \EE\left[ g(Z_1(h_1))| \bxi, Z_0=\delta_x \right] }{M_0h_1(x)}\right]\leq\EE \left[ \sup_{x\in \XX} \frac{ \sum_{k\geq 0} (k+1) \log^+(k+1)^{1+\varepsilon} \PP\left[F_{x,\xi_0}=k|\bxi \right] }{f_{x,\xi_0}}\right] .\]
	As a consequence, \eqref{eq:LlogL age} is indeed sufficient for \autoref{ass: LlogL}.
	Let us check now that \eqref{eq:LlogL age tail} also implies \autoref{ass: LlogL}.
	We perform an integration by parts in \eqref{eq: LogL avant IPP} and derive
	\begin{align*}
		\EE\left[ g(Z_1(h_1))| Z_0=\delta_x, \bxi \right]&\leq \sum_{k\geq 1} (g(k+1)-g(k)) \PP\left[F_{x,\xi_0}\geq k|\bxi \right].
	\end{align*}
	An elementary function study (see Appendix, Lemma \ref{lem: etude g}) allows to show that there exists some constant $A$ depending only on $\varepsilon>0$, such that $g(k+1)-g(k)\leq A \log(1+k)^{1+\varepsilon}$ for all $k\geq 1$. This implies 
	\begin{equation}\label{eq: LlogL apres IPP}
		\EE\left[ g(Z_1(h_1))| Z_0=\delta_x, \bxi \right]\leq A\sum_{k\geq 1} \log(1+k)^{1+\varepsilon} \PP\left[F_{x,\xi_0}\geq k|\bxi \right].
	\end{equation}
	Similarly, combining \eqref{eq: moyenne taille} and \eqref{eq: LlogL apres IPP}, one proves that 
	\[\EE \left[ \sup_{x\in \XX} \frac{ \EE\left[ g(Z_1(h_1))| \bxi, Z_0=\delta_x \right] }{M_0h_1(x)}\right]\leq A\EE \left[ \sup_{x\in \XX} \frac{  \sum_{k\geq 1} \log(1+k)^{1+\varepsilon} \PP\left[F_{x,\xi_0}\geq k|\bxi \right]}{f_{x,\xi_0}}\right], \]
	therefore \eqref{eq:LlogL age tail} indeed implies \autoref{ass: LlogL}.
\end{proof}
\color{black} \subsection{Asymptotic type distribution in the L-GWRE}
\label{subs: Leslie asymptotic type}
In the case of a L-GWRE, the condition \autoref{ass:L2} which is needed to apply Theorem \ref{thme:type_distrib} is satisfied as soon as the offpring distributions $(F_{x,e})_{x,e}$ satisfy the moment condition
$$\EE\left[\log^+\left( \sup_{x\in\XX} \EE[(F_{x,\xi_0}+1)^2]\right)\right]<\infty.$$
Since checking this last condition on a particular model is pretty straightforward, we do not discuss it further and focus now on the additionnal condition $$\EE[\log^+(\nu M_{1,N}(f))]<+\infty$$ which guarantees that the assertion \textit{ii)} of Theorem \ref{thme:type_distrib} holds.
Namely, we state and prove the following proposition.
\begin{proposition}\label{prop: type distribution LGWRE}
	Consider a L-GWRE such that Assumption \autoref{ass:leslie} holds. We assume moreover that $\EE[-\log(s_{x,\xi_0})]<+\infty$ for all $x\in\XX$. 
	Then for any non-negative and nonzero function $f$, there exists an deterministic integer $N\geq 0$ such that  $$\EE[\log^+(\delta_0 M_{1,N}(f))]<+\infty$$.\end{proposition}
When the conditions of this Proposition hold, as well as the moment conditions required for Theorem \ref{thme:type_distrib}, we can apply Theorem \ref{thme:type_distrib},\textit{ii)} to any function $f$, and in particular to $f=\mathds{1}_y$ for any type $y$. We obtain hence
\[\lim_{n\rightarrow \infty} \frac{Z_n(y)}{\EE[Z_n(y)|Z_0,\bxi]}=W \quad \PP\text{-a.s.}\]
as well as $\lim \left( Z_n(y)\right)^\frac{1}{n}={e^\lambda}$ almost surely on $\{W>0\}$, for any initial population $Z_0$ and any type $y\in\XX$.
\begin{proof}[Proof of Proposition \ref{prop: type distribution LGWRE}]
	Let $f$ be non-negative and non-zero. We can choose $y\in\XX$ such that $f(y)>0$ and write $f\geq f(y)\mathds{1}_y$. Thus, to prove $\EE[\log^-(\delta_0 M_{1,N} f)]<+\infty$ we only need to prove that $\EE[\log^{-}(M_{1,N}(0,y))]<+\infty$, where $M_{1,N}(0,y)$ refers to the entry on the $0$-th row and $y$-th column of the matrix $M_{0,N}$. We set $N=y+1$, and write $M_{1,y+1}(0,y)\geq M_{1}(0,1)\dots M_{y}(y-1,y)=\prod_{x=1}^{y} s_{x-1,\xi_x}$.
	Hence,
	$$\EE[\log^{-}(M_{1,N}(0,y))]\leq \sum_{x=1}^{y}\EE[\log^-(s_{x-1,\xi_x})]<+\infty,$$
	if $\EE[-\log(s_{x,\xi_0})]<+\infty$ for all $x$.
\end{proof}
\color{black}
\subsection{Extinction and explosion}\label{subs: Leslie extinction-explosion}
As explained in the introduction, we do not expect a dichotomy between the events $\ext$ and $\{W>0\}$ in general when $\XX$ is infinte. However, Proposition \ref{prop:extinction-explosion} states that this dichotomy holds when there exists a type $x_0$ such that $\limsup_{n\rightarrow \infty} Z_n(x_0)=+\infty$ a.s. on the survival event. In the next statement, we prove that under mild assumptions, this situation happens in the particular case of an L-GWRE.
\begin{prop}\label{prop:exp et ext LGW}
	Let $(Z_n)$ be a L-GWRE, satisfying the assumptions of Proposition \ref{prop:hyp_leslie_2}, as well as \autoref{ass: LlogL}. We suppose additionally that
	\begin{enumerate}[i)]
		\item \( \inf_{x,e} \PP[F_{x,e}=0]=a>0\)
		\item \(  s=\sup_{x,e} s_{x,e} <1.\)
	\end{enumerate}
	Then \eqref{eq:survie type} holds a.s. on survival. If additionally $\bxi$ is i.i.d then Proposition \ref{prop:extinction-explosion} applies and
	\[ \PP\left[ \nexp\triangle \ext\right]=0.\]
\end{prop}
Notice that under the stochastic monotonicity assumption \autoref{ass:leslie},\textit{ii)}, it holds $$a=\inf_{e\in\EE} \PP[F_{0,e}=0]\text{ and }s=\sup_{e\in\EE} s_{0,e}.$$
Our proof relies on two lemmas, in which we prove results of independent interest.
\begin{lem} \label{lem:ext-explo} Let $(Z_n)$ be a L-GWRE such that $\inf_{x,e} \PP[F_{x,e}=0]=a>0$ and $s=\sup_{x,e} s_{x,e} <1$,	
	then  \[ \PP_{\bxi}[\ext] + \PP_{\bxi}[ \lim_{n\rightarrow + \infty}\Vert Z_n \Vert =+\infty]=1 \quad \PP(d\bxi)\text{-a.s.} \]
\end{lem}
\begin{lem}\label{lem:expl-type0} 
	Let $(Z_n)$ be a L-GWRE such that $s=\sup_{x,e} s_{x,e} <1$,
	then it holds
	\[\left\{ \lim_{n\rightarrow\infty} \Vert Z_n \Vert =+\infty \right\}\subset \left\{\limsup_{n\rightarrow \infty} Z_n(0)=+\infty \right\}\] up to a $\PP$-negligeable event.\end{lem}
Once again, we postpone the proofs of the lemmas and focus first on the proof of Proposition \ref{prop:exp et ext LGW}, which is rather short.
\begin{proof}[Proof of Proposition \ref{prop:exp et ext LGW}.]
	Combining Lemmas \ref{lem:ext-explo} and \ref{lem:expl-type0}, we derive 
	\[\PP_{\bxi}[\ext] + \PP_{\bxi}[\limsup Z_n(0)=+\infty]\geq 1 \quad \PP(d\bxi)\text{-a.s}.\]
	Since the events $\ext$ and $\{\limsup Z_n(0)=+\infty\}$ are clearly incompatible, this proves
	\[\PP_{\bxi}[\ext] + \PP_{\bxi}[\limsup Z_n(0)=+\infty]= 1 \quad \PP(d\bxi)\text{-a.s},\]
	in other words, $\PP$-almost surely, conditionally on survival, it holds $\limsup_{n\rightarrow \infty} Z_n(0)=+\infty$. Thus Proposition \ref{prop:extinction-explosion} applies if $\bxi$ is i.i.d. 
\end{proof}
Let us move now to the proof of the lemmas.
\begin{proof}[Proof of Lemma \ref{lem:ext-explo}.]
	We adapt here arguments from \cite[Theorem 11.2]{harris_theory_1968}.
	We notice that
	\[\PP[Z_{1}=0| Z_0,\bxi]=\prod_{x\in\XX}\PP[(F_{x,\xi_0},S_{x,\xi_0})=0|\bxi]^{Z_0(x)}\geq (a s)^{\Vert Z_0 \Vert}.\]
	Let us note
	\begin{align*}
		R_K(\be,Z)&=\PP[ 0<\Vert Z_n \Vert \leq K \text{ infinitely often }|\bxi=\be, Z_0=Z ],\\
		\Omega_K&=\{Z\in \N | 0<\Vert Z \Vert\leq K\},\\
		\bar R_K&= \sup_{\substack{\be \in \E^{\ZZ_+} \\ Z\in \Omega_K}} R_K(\be,Z)\\
		\text{ and } T_K&=\inf\{n\geq 1, Z_n\in \Omega_K\}\in {\ZZ_+}\cup\{+\infty\},
	\end{align*}
	for each $K \geq 1$, $Z\in \N$ and $\be\in \E^{\ZZ_+}$.
	It clearly holds 
	\begin{align}R_K(\be,Z)&=\EE\left[R_K(\theta^{T_K}(\be),Z_{T_K})\mathds{1}_{{T_K}<\infty}|\bxi=\be, Z_0=Z\right] \nonumber\\&\leq \bar{R}_K \PP[T_K<\infty|\bxi=\be,Z_0=Z]. \label{eq:recurrence RK}\end{align}
	However,  $\{Z_1=0\}\subset \{T_K=\infty\}$, thus \[\PP[T_K<\infty|\bxi=\be,Z_0=Z]\leq 1-\PP[Z_{1}=0|\bxi=\be,Z_0=Z]\leq 1-(a s)^{\Vert Z \Vert }.\]
	In particular, if $Z\in \Omega_K$ then $\Vert Z \Vert \leq K$ thus
	\[\PP[T_K<\infty|\bxi=\be,Z_0=Z]\leq 1-(a s)^K.\]
	Plugging this into \eqref{eq:recurrence RK}, and taking a supremum over $Z\in\Omega_K$, we get
	\[\bar{R}_K \leq \bar{R}_K \sup_{\substack{\be \in \E^{\ZZ_+} \\ Z\in \Omega_K}}\PP[T_K<\infty|\bxi=\be,Z_0=Z] \leq \bar{R}_K \left(1-(a s)^K\right),\]
	where $1-(a s)^K<1$. Thus $\bar{R}_K=0$ for any $K\geq 1$. Hence for any $Z\in \N$ and $\be \in \E^{\ZZ_+}$ it holds
	\[  \PP\left[0<\liminf \Vert Z_n \Vert  <\infty | \bxi=\be, Z_0=Z\right]=\lim_{K\rightarrow \infty} R_K(\be,Z)=0,\]
	which concludes.
\end{proof}
\begin{proof}[Proof of Lemma \ref{lem:expl-type0}.]
	By Definition \ref{eq:MGWRE}, the L-GWRE $(Z_n)$ satisfies the recurrence relation
	\begin{align*} Z_{n+1}(x+1)=\sum_{k=1}^{Z_n(x)} S_{x,\xi_n}^{k,n} \text{ and }
		Z_{n+1}(0)=\sum_{x\in\XX} F_{x,\xi_n}^{k,n},\end{align*}
	the $(S_{x,e}^{k,n} : e\in \E, k,n,x\geq 0)$ are mutually independent Bernoulli variables such that $\EE[S_{x,e}^{k,n}]=s_{x,e}\leq s$. Thus, we may consider a coupling of the array $(S_{x,e}^{k,n})_{k,n,x,e}$ with an array $(\tilde{S}_{x,e}^{k,n})_{x,k,n,e}$ of mutually independent Bernoulli variables of mean $s$ in such a way that almost surely, for any $k,n\in \ZZ_+,x\in\XX,e\in\E$
	\[ S_{x,e}^{k,n} \leq \tilde{S}_{x,e}^{k,n}.\]
	For each integer $K$, we define a process $(\tilde{Z}^K_n)$ with values in $\N$ by setting
	$\tilde{Z}_0=Z_0$ and for any $n\geq 0$,
	\[ \tilde{Z}^K_{n+1}(x+1)=\sum_{i=1}^{\tilde{Z}_{n}^K(x)} \tilde{S}_{x,\xi_n}^{i,n} \text{ for }x\geq 0,\text{ and }
	\tilde{Z}^K_{n+1}(0)= K.\]
	Recalling the notation $\Omega_K=\{ \sup_{n\geq 0} Z_n(0) \leq K \},$ a mere induction allows to show that for all $n\in{\ZZ_+}$, we have $Z_n \mathds{1}_{\Omega_K}\leq \tilde{Z}_n^K.$
	Moreover, it holds
	\[ \Vert \tilde{Z}^K_{n+1} \Vert=K+\sum_{x\in\XX} \sum_{i=1}^{Z_n(x)} S_{x,\xi_n}^{i,n},\]
	and conditionally on $\bxi$ and $\tilde{Z}_n^K$
	\[\Vert \tilde{Z}^K_{n+1} \Vert \overset{d}{=} K+\sum_{i=1}^{\Vert \tilde{Z}_n^K\Vert} S^{i,n}, \]
	where the $S^{i,n}$ are independent Bernoulli variables of parameter $s$.
	Therefore, the process $(\Vert Z^K_{n}\Vert)_{n\geq 0} $ is distributed as a Galton-Watson process with immigration, where the offspring distribution is a Bernoulli law of parameter $s$ and the immigration is deterministic with value $K$. This process describes a monotype particle system where each particle present at time $n$ survives independently from the others with probability $s$ (and dies with probability $1-s$), and $K$ new particles arrive at each time step. As a consequence, noting $S'^n_{i,k}$ the Bernoulli variable describing the survival from time $n$ to time $n+1$ of the $i$-th particle arrived at time $k$, it is clear that the law of the process $(\tilde{Z}^K_n)_{0\leq n}$ satisfies 
	\[(\tilde{Z}^K_n)_{n\geq 0} \overset{d}{=} \left(\sum_{k=1}^{n} \sum_{i=1}^K S'^k_{i,k}\dots S'^{n-1}_{i,k} + \sum_{i=1}^{\Vert \tilde{Z}_0^K\Vert } S'^{0}_{i,0}\dots S'^{n-1}_{i,0}.\right)_{n\geq 0},\]
	where the $(S'^n_{i,k})_{i,k,n\geq 0}$ are mutually independent. 
	At a fixed time $n$, it holds therefore
	\begin{align*} \tilde{Z}^K_n\overset{d}{=}\sum_{k=1}^{n} \sum_{i=1}^K S'^k_{i,k}\dots S'^{n-1}_{i,k} + \sum_{i=1}^{\Vert \tilde{Z}_0^K\Vert } S'^{0}_{i,0}\dots S'^{n-1}_{i,0}\overset{d}{=} \sum_{k=0}^{n-1} \sum_{i=1}^K S'^1_{i,k}\dots S'^{k}_{i,k} + \sum_{i=1}^{\Vert \tilde{Z}_0^K\Vert } S'^{0}_{i,0}\dots S'^{n-1}_{i,0}.\end{align*}
	And as $n\rightarrow \infty,$ it holds $\PP$-almost surely
	\[\sum_{k=0}^{n-1} \sum_{i=1}^K S'^1_{i,k}\dots S'^{k}_{i,k} + \sum_{i=1}^{\Vert \tilde{Z}_0^K\Vert } S'^{0}_{i,0}\dots S'^{n-1}_{i,0} \underset{n\rightarrow \infty}{\longrightarrow} \sum_{k=0}^{+\infty} \sum_{i=1}^K S'^1_{i,k}\dots S'^{k}_{i,k}=:S^K.\]
	Therefore, 
	\[ \Vert \tilde{Z}^K_n \Vert \overset{d}{\underset{n\rightarrow \infty}{\longrightarrow}} S^K.\]
	Additionally, note that $\EE[S^K]=\sum_{i=1}^K \sum_{k=0}^{+\infty} \EE[S'^1_{i,k}\dots S'^{k}_{i,k}]\leq K/(1-s)<\infty$ thus $S^K$ is $\PP$-a.s. finite.
	
	In particular, since $\Vert \tilde{Z}_n^K\Vert$ and $S^K$ are integer-valued, this convergence implies 
	\begin{align*}
		\PP\left[\liminf_{n\rightarrow\infty}\Vert \tilde{Z}^K_n \Vert =+ \infty \right] &= \PP\left[\bigcap_{A\in \ZZ^*_+} \bigcup_{N\geq 0} \bigcap_{n\geq N}  \left\{ \Vert \tilde{Z}^K_n \Vert >A \right\}\right] \\
		&=\liminf_{A\rightarrow \infty} \limsup_{N\rightarrow \infty} \PP\left[\bigcap_{n\geq N} \left\{  \Vert \tilde{Z}^K_n \Vert>A\right\}\right] \\
		&\leq \liminf_{A\rightarrow \infty} \limsup_{N\rightarrow \infty} \PP\left[ \Vert \tilde{Z}^K_N \Vert >A\right] \\
		&\leq \liminf_{A\rightarrow \infty} \PP\left[S^K>A \right] \\
		&=0.\end{align*}
	Since once again $Z_n\mathds{1}^{\Omega^K}\leq \tilde{Z}_n^K$, it holds therefore 
	\[\PP[\liminf_{n\rightarrow +\infty} \mathds{1}_{\Omega_K}\Vert Z_n\Vert  =+ \infty]=0.\]
	Therefore, $\PP$-almost surely, either 
	\[\forall K\geq 1, \,\sup_{n\geq 1} Z_n(0)>K, \text{ i.e. } \limsup_{n\rightarrow\infty} Z_n(0)=+\infty\]
	or
	\[\liminf_{n\rightarrow\infty} \Vert Z_n \Vert <+\infty.\]
	Therefore, up to a $ \PP$-negligible event, it holds
	\[ \{ \lim \Vert Z_n \Vert=+\infty \} = \{ \liminf \Vert Z_n \Vert = +\infty\} \subset \{ \limsup Z_n(0)=+\infty\}. \]
\end{proof}

\section{Appendix}
\begin{lemap}\label{lem: etude u}
	Let $u:(a,s)\in (1+\infty)^2\mapsto\sum_{k\geq 1} \left(\log k\right)^sk^{-a}$. 
	For all $a,s>1$, it holds $$u(a,s)=\frac{\Gamma(s+1)}{(a-1)^{s+1}}+2 \left(\frac{s}{e}\right)^s.$$
\end{lemap}
\begin{proof}
	We note $\psi_{s,a}(t)=\log(t)t^{-a/s}$ and notice that $u(a,s)=\sum_{k\geq 1} \psi_{s,a}(k)^s.$
	The derivative of $\psi_{s,a}$ is 
	\[\psi_{s,a}'(t)=t^{-\frac{a}{s}-1} \left(1-\frac{a}{s}\log(t)\right),\]
	Thus $\psi_{s,a}$ and $\psi_{s,a}^s$ are increasing on $(1,e^{s/a}]$ and decreasing on $(e^{s/a},+\infty)$.
	Therefore, we write
	\[u(s,a)\leq \sum_{1 \leq k \leq e^{s/a}-1} \psi_{s,a}(k)^s+2\psi_{s,a}(e^{s/a})^s+\sum_{e^{s/a}+1\leq k} \psi_{s,a}(k)^s.\]
	Where
	\[\sum_{1 \leq k \leq e^{s/a}-1} \psi_{s,a}(k)^s\leq \int_{1}^{e^{s/a}} \psi_{s,a}(t)^sdt\]
	and
	\[\sum_{1+e^{s/a}\leq k} \psi_{s,a}(k)^s\leq \int_{e^{s/a}}^{+\infty} \psi_{s,a}(t)^s dt.\]
	This yields 
	\[u(s,a)\leq \int_1^{+\infty} \left(\log t\right)^st^{-a} dt + 2 \left(\frac{s}{a}\right)^s \exp(s/a)^{-a}.\]
	On the one hand, setting $z=(a-1)\log(t)$, one gets
	\begin{align*}
		\int_1^{+\infty} \log(t)^st^{-a} dt& =\int_{1}^{+\infty} \log(t)^s e^{-(a-1)\log(t)}\frac{dt}{t} \\
		&=\frac{1}{(a-1)^{s+1}}\int_{0}^{+\infty}z^se^{-z}dz \\
		&=\frac{\Gamma(s+1)}{(a-1)^{s+1}}
	\end{align*}
	On the other hand 
	\[ \left(\frac{s}{a}\right)^s \exp(s/a)^{-a}=\left(\frac{s}{ae}\right)^s\leq \left(\frac{s}{e}\right)^s\]
	since $a>1$.
	Therefore 
	\[u(s,a)\leq \frac{\Gamma(s+1)}{(a-1)^{s+1}}+2 \left(\frac{s}{e}\right)^s\]
	which is indeed the desired result.
\end{proof}

\begin{lemap}\label{lem: etude g}
	Let $\varepsilon>0$ and $g(t)=t\left(\log t\right)^{1+\varepsilon}$ for $t\geq 1$.
	Then there exists $A>0$ such that $g(k+1)-g(k)\leq A \left(\log(k+1)\right)^{1+\varepsilon}$ for any $k\geq 1$.
\end{lemap}
\begin{proof}
	We compute the derivative
	\[g'(t)=\left(\log t\right)^{1+\varepsilon}+(1+\varepsilon)\left(\log t\right)^{\varepsilon}, \quad t\geq 1.\]
	By the mean value theorem, for $t\geq 1$, it holds therefore 
	\begin{align*}g(t+1)-g(t)& \leq \sup_{t\leq s \leq t+1} g'(s) \\
		&\leq  g'(t+1)
		\\& \leq \left(\log(t+1)\right)^{1+\varepsilon}+(1+\varepsilon)\left(\log(t+1)\right)^{\varepsilon} 
		\\& \leq \left(1+\frac{1+\varepsilon}{\log(2)}\right) \left(\log(1+t)\right)^{1+\varepsilon}.\end{align*}
	This holds therefore in particular for $t=k\in \ZZ_+.$
\end{proof}

\section*{Acknowledgements}
I would like to warmly thank my PhD-supervisors Vincent Bansaye and Marc Peigné for their continous guidance and support as well as their numerous feedbacks on this manuscript, which greatly helped improving it.

\printbibliography
\end{document}